\documentclass[reqno,11pt]{amsart}
\usepackage{amssymb,latexsym}

\newtheorem{thm}{Theorem}
\newtheorem{lemma}{Lemma}
\newtheorem{cor}{Corollary}
\newtheorem{prop}{Proposition}

\theoremstyle{definition}

\theoremstyle{remark}

\makeatletter
\@addtoreset{figure}{section}
\def\thefigure{\thesection.\@arabic\c@figure}
\def\fps@figure{h, t}
\@addtoreset{table}{bsection}
\def\thetable{\thesection.\@arabic\c@table}
\def\fps@table{h, t}
\@addtoreset{equation}{section}

\makeatother

\allowdisplaybreaks
\def\intprod{\mathbin{\hbox to 6pt{%
                 \vrule height0.4pt width5pt depth0pt
                 \kern-.4pt
                 \vrule height6pt width0.4pt depth0pt\hss}}}

\begin{document}

\title{On Incompressible Averaged Lagrangian Hydrodynamics}

\author{Steve Shkoller}

\address{Department of Mathematics\\
University of California \\
Davis, CA 95616}
\email{shkoller@math.ucdavis.edu}


\date{May 29, 1999; current version September 28, 1999}

\maketitle

\tableofcontents

\section{Introduction}
\label{Intro}
The Euler equations of incompressible hydrodynamics and their viscous
counterparts, the Navier-Stokes equations, arise as the lowest order models in
a hierarchy of fluids models known as differential-type fluids or
Rivlin-Ericksen fluids developed in \cite{RE} forty-five years ago.  
A first order correction to these
equations is provided by the {\it second-grade fluids} model  \cite{NT},
written in Euclidean space  ${\mathbb R}^n$ as the following system:
\begin{equation}\label{Euclidean}
\begin{array}{c}
\partial_t (1-\alpha^2 \triangle) u - \nu\triangle u +
\text{curl}(1-\alpha^2\triangle)u \times u = -\text{grad } p, \\
u(0) = u_0, \ \ \ \text{div }u=0,
\end{array}
\end{equation}
where $u(t,x)$ is the velocity vector field, $p(t,x)$ is the scalar pressure,
and $\alpha, \nu$ are positive constants.  Formally, as $\alpha \rightarrow 0$,
the Navier-Stokes equations are recovered.

Although the system (\ref{Euclidean}) has been widely studied in the
mathematics literature (see \cite{CO}, \cite{CG}, \cite{GGS}, \cite{DF},
\cite{RE}, \cite{NT}, and the many references therein), its
complicated nonlinearity, mixed temporal-spatial differential
operators, and incompressibility constraint have caused difficulties for
traditional analytic techniques, and the following fundamental problems 
have remained
open: Lagrangian boundary conditions, local existence and uniqueness of 
smooth-in-time solutions for
the inviscid ($\nu=0$) problem in dimension three, viscosity independent
time intervals of existence, and the regular limit of zero viscosity 
in fluid domains with boundary.   In this paper, we shall state three
main theorems which solve these problems.  Our method is geometric, and
relies heavily on properties of certain nonlinear operators between sections
of infinite-dimensional Hilbert bundles over new subgroups of the 
volume-preserving diffeomorphism group.

We have been motivated to study the equation (\ref{Euclidean}) because
of the remarkable fact that the recently developed {\it averaged
Euler equation} or Euler-$\alpha$ equation (see, for example,  
\cite{HMR1}, \cite{HMR2}, \cite{S}, 
\cite{MRS}, \cite{MS}), which was introduced as an 
LES\footnote{LES stands for Large Eddy Simulation and constitutes
a class of models that average over the small scales of the fluid which
cannot be resolved computationally.}-type mathematical model for incompressible
fluid flow, is mathematically identical
to the second-grade fluids equation with zero viscosity; with viscosity
present (\ref{Euclidean}) is known as the averaged Navier-Stokes equation.  
This coincidence
is completely surprising because the parameter $\alpha$ in the Rivlin-Ericksen
hierarchy represents a material parameter measuring the elastic response
of the fluid, while this parameter in the averaged Euler-formulation
denotes a spatial length scale.  In yet another striking coincidence,
the equation (\ref{Euclidean}) exactly coincides with the vortex blob numerical
algorithm introduced by Chorin in \cite{C}  with a particular choice of 
smoothing function.  Irregardless of the context in which (\ref{avg_euler})
is considered, solutions of this equation with sufficiently small
$\alpha$  are able to qualitatively reproduce
the behaviour of the large scale flow (spatial scales larger than $\alpha$)
for high-Reynolds number incompressible fluids, while filtering or averaging
over the small-scales \cite{NS} \cite{CHMZ}.  As a result this model is
better suited for numerical simulations of complex fluid flow and turbulence
\cite{FHT2}.  We study the analytic and geometric properties of this fluid
motion.

We first generalize the system of equations (\ref{Euclidean}) from
${\mathbb R}^n$  to the setting of a $C^\infty$  compact oriented
$n$-dimensional Riemannian  manifold with $C^\infty$ boundary, $(M,g)$.
Letting $\nabla$ denote the Levi-Civita covariant derivative, 
(\ref{Euclidean}) becomes
\begin{equation}\label{avg_euler}
\begin{array}{c}
\partial_t(1-\alpha^2\triangle_r)u - \nu \triangle_r u + 
\nabla_u(1-\alpha^2\triangle_r)u
-\alpha^2 (\nabla u)^t \cdot \triangle_r u = -\text{\text{grad }}p,\\
\text{\text{div }} u =0, \ \ u(0)=u_0,\\
\alpha >0, \ \ \ \triangle_r = -(d \delta + \delta d) + 2\text{Ric},
\end{array}
\end{equation}
together with one of the following three boundary conditions:
\begin{itemize}
\item[(a)] Dirichlet or no-slip: $u=0$ on ${\partial M}$,
\item[(b)] Neumann or free-slip: 
$g(u,n)=0$ and $(\nabla_nu)^{\rm tan} + S_n(u)=0$
on ${\partial M}$,
\item[(c)] mixed: $u=0$ on $\Gamma_1$, and
$g(u,n)=0$, $(\nabla_nu)^{\rm tan} + S_n(u)=0$ on $\Gamma_2$,
where ${\partial M} = \Gamma_1\cup \Gamma_2$,
$\overline{\Gamma}_1 = {\partial M}/\Gamma_2$, and the sets $\Gamma_1$,
$\Gamma_2$ are disjoint.
\end{itemize}
On a Riemannian manifold,
there is always more than one choice for the correct ``Laplacian'' on
vector fields or $1$-forms.  Our Laplacian $\triangle_r$ is the
operator ${\mathcal L}= -2\text{Def}^*\text{Def}$ acting on 
divergence-free vector fields (or coexact $1$-forms), where 
the (rate of) deformation tensor is given by
$$\text{Def}(u)= {\frac{1}{2}}(\nabla u + \nabla u ^t) = {\frac{1}{2}}
\pounds _u g,
\ \ \ \pounds = \text{ Lie derivative},$$
and $\text{Def}^*$ is the $L^2$ formal adjoint of $\text{Def}$.
Other possible choices are the Hodge Laplacian $-\triangle=
(d \delta + \delta  d)$ or the rough Laplacian $-\text{Tr}\nabla \nabla$,
but the boundary conditions (a)-(c) insist upon our choice ${\mathcal L}$.
Note that
$$(\nabla_nu)^{\rm tan} + S_n(u)= [\text{Def}(u)\cdot n]^{\rm tan}$$
for vector fields $u$ which are tangential to ${\partial M}$.

\noindent{\bf Further Notation.}
For each $x \in {\partial M}$, the $g$-orthogonal bundle splitting
$T_xM = T_x{\partial M} \oplus N_x$ induces the Whitney sum
$$ TM|_{\partial M} = T\partial M \oplus_g N,$$
where $N$ is the normal bundle, $N=\cup_{x\in {\partial M}}N_x \downarrow
{\partial M}$.

Letting $\pi:E \rightarrow M$ be a vector bundle over $M$ (or over
${\partial M}$), we denote the $H^s$ sections of $E$ by $H^s(E)$ and
for all $\eta \in {\mathcal D}^s$, we set
$H^s_\eta(E):= \{ U \in H^s(M,E) \mid \pi \circ U = \eta\}$.

For any vector bundle ${\mathcal E}$ over a base manifold ${\mathcal M}$,
we shall often make use of the notation 
${\mathcal E}_m \downarrow {\mathcal M}$ to denote ${\mathcal E}$, where 
${\mathcal E}_m$ is the fiber over $m \in {\mathcal M}$.

We use $R$ to denote the Riemannian curvature operator of $\nabla$.
The Ricci curvature as a bilinear form is given by
$$\text{Ric}(x,y)= \text{Tr }g(R(x,\cdot )\cdot,y ),$$
with the associated linear operator Ric$:T_xM \rightarrow T_xM$ given by
$g(\text{Ric}(x),y) = \text{Ric}(x,y)$.

If $\eta \in H^s(M,M)$, then the tangent mapping $T\eta$ is in $H^{s-1}(M,
T^*M \otimes \eta^*(TM))$.  If $w\in T_xM$, then in a local
chart, $T\eta(x) \cdot w = (\eta(x), D\eta(x) \cdot w)$ where $D$ is the matrix
of partial derivatives of $\eta$ with respect to the coordinate chart.

We shall use the symbol $\pounds$ to denote the Lie derivative, $d$ for
the exterior derivative on $\Lambda^k(M)$, the differential $k$-forms on
$M$, and $\delta$ for its
formal adjoint with respect to the $L^2$ pairing.  For a vector field $u$
on $M$, $\nabla u^t$ shall denote the transpose of $\nabla u$ with respect to
$g$.  

The Hodge Laplacian on differential $k$-forms is $\triangle =
-(d \delta + \delta d)$, and 
$$\triangle_r = \triangle + 2\text{Ric}.$$

When we wish to explicitly convert between vector fields and $1$-forms,
we shall use the
musical maps $\flat : TM \rightarrow T^*M$ and $\sharp:T^*M \rightarrow TM$;
for example, if $u$ is a vector field on $M$, then $u^\flat$ is the
associated $1$-form.

\section{Main Results}
We prove the existence of smooth-in-time classical solutions to
(\ref{avg_euler}) by transforming the Eulerian equation given above into
a Cauchy problem for the Lagrangian flow map on any one of three subgroups of 
${\mathcal D}_\mu^s$,
the topological group consisting of Hilbert $H^s$-class volume-preserving
diffeomorphism of $M$ with $H^s$ inverses.  Our
first theorem proves the existence of these subgroups.

\begin{thm}\label{Thm1}
We take $s>{\frac{n}{2}}+1$.  For
$\eta \in {\mathcal D}_\mu^s$, let $T\eta$ denote its tangent map,
i.e., the Frechet derivative of $\eta$ thought of as bundle map. 
Let $n$ denote the outward-pointing normal field along the boundary
 $\partial M$, and
let $S_n:T{\partial M} \rightarrow T{\partial M}$ denote the Weingarten
map so that
$$g\left( S_n(u), v \right) = II_n(u,v)=-g(\nabla_un,v), \ \
u,v \in H^{s-{\frac{3}{2}}}(T{\partial M}),$$
where $II_n$ is the second fundamental form of ${\partial M} \subset M$.
Define the sets
\begin{align*}
{\mathcal D}^s_{\mu,N} & = \{ \eta \in {\mathcal D}_\mu^s \mid
T\eta|_{\partial M} \cdot n \in H^{s-3/2}_\eta(N),
\text{ \rm for all } n \in H^{s-1/2}(N)\},\\
{\mathcal D}^s_{\mu,D} & =\{ \eta \in {\mathcal D}_\mu^s\ | \
\eta|_{\partial M} = e\},
\end{align*}
and
\begin{align*}
{\mathcal D}^s_{\mu,mix} & =\{ \eta \in {\mathcal D}_\mu^s\ | \ \eta
\text{ \rm leaves } \Gamma_i \text{ \rm invariant}, \ \ \eta|_{\Gamma_1} = e,
 \\& \qquad 
T\eta|_{\Gamma_2} \cdot n \in H^{s-3/2}(N|\Gamma_2), \text{ \rm for all }
n\in H^{s-1/2}(N|\Gamma_2) \},
\end{align*}
where
we suppose that $M,{\partial M}$ are $C^\infty$, that $\Gamma_1$ and $\Gamma_2$
are two {\it disjoint} subsets of ${\partial M}$ such that if $m_0 \in
\Gamma_i$ $(i=1,2)$, a local chart $U$ (in $\overline M$) about $m_0$ can be
chosen so that $\overline U \cap {\partial M} \subset \Gamma_i$; furthermore,
we assume that $\overline \Gamma_1 = {\partial M}/\Gamma_2$ and that
${\partial M} = \Gamma_1 \cup \Gamma_2$.

Then ${\mathcal D}_{\mu,D}^s$, ${\mathcal D}_{\mu,N}^s$, and 
${\mathcal D}_{\mu,mix}^s$ are all $C^\infty$ subgroups of ${\mathcal D}_\mu^s$,
and the tangent space at the identity of these groups is given by
\begin{align*}
T_e {\mathcal D}^s_{\mu,N} = \{ u \in &T_e {\mathcal D}_\mu^s \mid
(\nabla_n u|_{\partial M})^{\rm tan} + S_n(u)=0\in
H^{s-3/2}(T{\partial M}) \\
&\text{\rm for all } n \in H^{s-1/2}(N) \},\\
T_e {\mathcal D}^S_{\mu,D} = \{ u \in &T_e {\mathcal D}_\mu^s \mid
u|_{\partial M}=0 \}, 
\end{align*}
and
\begin{align*}
T_e {\mathcal D}^S_{\mu,mix} &= \{ u \in T_e {\mathcal D}_\mu^s \mid
(\nabla_n u|_{\partial M})^{\rm tan} + S_n(u)=0\in
H^{s-{\frac{3}{2}}}(T\Gamma_2) \\
&\qquad \text{\rm for all } n \in
H^{s-{\frac{1}{2}}}(N|\Gamma_2)\text{  \rm and  } u|_{\Gamma_1}=0 \}.
\end{align*}

We also form the corresponding sets ${\mathcal D}_{N}^s$,
${\mathcal D}_{D}^s$, and ${\mathcal D}_{mix}^s$ which do not have the
volume-preserving constraint imposed. These sets are $C^\infty$ subgroups
of the full diffeomorphism group ${\mathcal D}^s$, and have the
analogous tangent spaces at the identity without the divergence-free
constraint.
\end{thm}

We call the groups ${\mathcal D}_{\mu,D}^s$, ${\mathcal D}_{\mu,N}^s$,
and ${\mathcal D}_{\mu,mix}^s$, the Dirichlet, Neumann, and mixed 
volume-preserving diffeomorphism groups.  Theorem \ref{Thm1} allows us
to do smooth calculus on these spaces.  We can thus transform the rather
complicated evolution equation (\ref{avg_euler}) to a simpler Cauchy
problem for the Lagrangian flow on these spaces.  In this article, we
shall prove results for the case of the group ${\mathcal D}_{\mu,D}^s$, as
the no-slip conditions have been of most interest in the literature.
We are able to prove the following result.\footnote{By setting 
${\partial M}=\emptyset$, all of our results hold for boundaryless 
compact oriented $n$ dimensional Riemannian manifolds.}

\begin{thm}\label{Thm2}
Set $s>(n/2)+1$,  and
let $\langle \cdot, \cdot \rangle$ denote the right invariant metric
on ${\mathcal D}_{\mu,D}^s$ given at the identity by
$$\langle X,Y \rangle_e =  (X,Y)_{L^2} +
{\frac{\alpha^2}{2}} (\pounds_X g, \pounds_Yg)_{L^2}.$$
For $u_0 \in T_e {\mathcal D}_{\mu,D}^s $,  there exist 
intervals $I=(-T,T)$ and $\bar I=[0,T)$, depending on $|u_0|_s$, and
a unique geodesic $\dot \eta$ of $\langle \cdot, \cdot \rangle$
with initial data
$\eta(0) = e$ and $\dot \eta(0)=u_0$ such that
$$\dot \eta \text{ is in } C^\infty( I, T{\mathcal D}^s_{\mu,D})$$
and has $C^\infty$  dependence on the initial velocity $u_0$.  

The geodesic $\eta$ is the Lagrangian flow of the time-dependent vector
field $u(t,x)$ given by
$$ \partial_t\eta(t,x) = u(t,\eta(t,x)),$$
and 
$$u \in  C^0( I, {\mathcal V}^s_\mu) \cap C^1(I,{\mathcal V}^{s-1}_\mu)$$
uniquely solves (\ref{avg_euler}) with $\nu=0$, and depends continuously
on $u_0$.

Furthermore, if
for $r\ge 1$, we set ${\mathcal V}^r_\mu =\{ u \in H^s(TM) \cap H^1_0(TM)\ | \
\text{\rm \text{div }}u=0\}$,
then 
$$ u^\nu \in C^0( \bar I, {\mathcal V}^s_\mu) \cap 
C^1(\bar I,{\mathcal V}^{s-1}_\mu)$$
is the unique solution of (\ref{avg_euler}) for $\nu >0$, depends continuously
on $u_0$, and $T$ is independent of the viscosity $\nu$.
\end{thm}

This is the first analytic result for solutions of (\ref{avg_euler}) with
or without viscosity that gives smooth evolution curves and smooth
dependence on initial data (c.f., \cite{CO},\cite{CG} and \cite{GGS}).
\begin{cor}\label{Cor1}
In the case that the dimension of $M$ is equal to two, $T=\infty$ and does
not depend on $\nu$.
\end{cor}
\begin{proof}
When $u$ is thought of as a $1$-form field on $M$, (\ref{avg_euler}) may be
reexpressed as
\begin{equation}\nonumber
\partial_t(1-\alpha^2 \triangle_r)u + \pounds_u(1-\alpha^2\triangle_r)u
=-d\hat p.
\end{equation}
Taking the exterior derivative of this equation and setting 
$q(t,x) = d (1-\alpha^2 \triangle_r) u(t,x)$
yields the vorticity equation
\begin{equation}\nonumber
\partial_t q +\pounds_u q =0.
\end{equation}
On two-dimensional manifolds,  we may identify the $2$-form $q$ with its
scalar density, in which case the above equation takes the particularly
simple form
\begin{equation}\label{vortex}
\begin{array}{c}
\partial_t q(t,x) +g(u(t,x),\operatorname{grad} q(t,x)) =0,\\
\operatorname{div}u=0,  \ \ u=0 \text{ on } {\partial M},\\
q(0)=d(1-\alpha^2\triangle_r)u_0,
\end{array}
\end{equation}
with the corresponding weak form
\begin{equation}\label{weak}
\int_{\mathbb R} \int_M \left(
q(t,x) \cdot \partial_t \phi(t,x) + g(q(t,x)u(t,x), \text{grad }\phi(t,x))
\right) \mu(x) dt =0
\end{equation}
for all $\phi \in C^\infty({\mathbb R} \times M)$.
Equation (\ref{vortex})  is equivalent
to the pointwise conservation of vorticity along the Lagrangian trajectory
\begin{equation}\label{pw}
q_t \circ \eta_t = q_0;
\end{equation}
 this is, of course, just the coadjoint action of 
${\mathcal D}_{\mu,D}^s$ acting on $T_e{\mathcal D}_{\mu,D}^s$ by right
composition.

Theorem \ref{Thm2} gives a time interval $[-T,T]$ of existence of solutions 
to (\ref{vortex}) for $u_0 \in {\mathcal V}^3$.  This, in turn, gives
the existence of a weak solution $q$ in $L^2(TM)$.
The conservation law (\ref{pw}) yields the Casimirs
$$ \int_M q^n \mu, \ \ n={\mathbb N}.$$
Thus, we have that $\|q(t,\cdot)\|_{L^2(M)}$ is a conserved quantity, and by 
standard elliptic estimates we have that $\|u(t,\cdot)\|_{H^3(TM)} 
\le K$ for all $t$.  Thus, the $H^3$-norm of $u(t,\cdot)$ does not blow up, 
so $T=\infty$.

The usual bootstrap argument now yields the result for $u_0 \in {\mathcal V}^s$,
$s>3$.
\end{proof}

This corollary is certainly not sharp, but rather demonstrates the ease with
which one obtains global-in-time solutions for (\ref{avg_euler}) in 2D for 
smooth enough initial data;  the proof of global existence for 2D Euler is
much more difficult, because one must rely on $L^\infty$ control
of vorticity and very careful $L^p$ estimates relying on quasi-Lipschitz
inequalities.

For (\ref{weak}), unique global solutions exist even for point-vortex 
initial data in the space of Radon measures
\cite{OS}.  Thus, one can solve a point-vortex ODE and generate
a {\it unique} global PDE solution.
This is not known to be the case for the point-vortex ODE associated
with the 2D Euler equation, for which the least regular initial data that gives
weak solutions is a vortex sheet (see \cite{Del}).

As a consequence of $T$ being independent of the viscosity $\nu$ for solutions
to (\ref{avg_euler}), we immediately obtain the following:
\begin{cor}\label{Cor2}
For $s>(n/2)+1$, solutions $u^\nu$ of (\ref{avg_euler}) converge regularly to 
the inviscid
solutions $u$ as $\nu \rightarrow 0$. Furthermore, letting 
$u^\nu = \partial_t \eta_\nu \circ \eta_\nu^{-1}$,  the viscous Lagrangian flow
$\eta_\nu$ converges regularly in the $H^s$ topology to the geodesic flow 
$\eta$ of the right invariant metric $\langle \cdot,\cdot \rangle$.
\end{cor}
This result states that we can generate smooth-in-time solutions to
(\ref{avg_euler}) with $\nu=0$ by obtaining a sequence of viscous solutions
with $\nu$ tending to zero.  Locally Lipschitz solutions were generated
in \cite{MRS}, so this result provides a significant improvement.
Thus, our result proves that the flow of the averaged Navier-Stokes equation
converges to the flow of the averaged Euler equation even in the presence
of boundaries.  This is in agreement with the scaling arguments of 
Barenblatt-Chorin (see, for example, the second paragraph of \cite{BC}).

We remark that traditional techniques, employed in \cite{CO}, \cite{CG}, 
have crucially relied
on viscosity to obtain existence of classical solutions. The critical
estimates in those papers have $1/\nu$ bounds, which prevent a limit of zero
viscosity result.  

The viscous term in (\ref{avg_euler}) is given by $-\nu \triangle_r u$,
and is derived from a rather deep constitutive theory for simple materials
\cite{NT}.  It is possible, however, to study this system on domains
without boundary, with stronger forms of viscosity.  For example, on
the three-dimensional torus, the article \cite{FHT} uses the 
dissipative term $-\nu\triangle(1-\alpha^2\triangle)u$ instead, which is
strong enough to guarantee global-in-time existence and uniqueness.
Following the product formula approach
developed in \cite{EM}, we can prove a regular limit of zero viscosity for
this type of dissipation as well.  For the following, which is Theorem 13.1 in
\cite{EM},  we assume that $M$ has no boundary.

\begin{prop}\label{propEM}
Let ${\mathcal B}: T\mathcal{D}_{\mu}^s\rightarrow T^2\mathcal{D}_{\mu}^s$
be the $C^\infty$ geodesic spray of the metric $\langle\cdot,\cdot\rangle_1$.
For each $s>(n/2)+1$ and $\sigma \ge 2$, 
let ${\mathcal T}:T_e \mathcal{D}_{\mu}^{s}
\rightarrow T_e \mathcal{D}_{\mu}^{s-\sigma}$ be
a bounded linear map that generates a strongly-continuous
semi-group $F_t:T_e\mathcal{D}_{\mu}^s\rightarrow T_e\mathcal{D}_{\mu}^s$,
$t\ge 0$, and satisfies $\|F_t\|_s \le e^{\beta t}$ for some $\beta > 0$
and some $s$.
Extend $F_t$ to $T\mathcal{D}_{\mu}^s$ by
$$ \tilde{F}_t(X_\eta) = TR_\eta \cdot F_t \cdot TR_{\eta^{-1}} (X_\eta)$$
for $X_\eta \in T_\eta \mathcal{D}_{\mu}^s$, and let
$\tilde{{\mathcal T}}$ be the vector
field $\tilde{{\mathcal T}}:T\mathcal{D}_{\mu}^s \rightarrow
T^2\mathcal{D}_{\mu,0}^{s-\sigma}$ associated to the flow $\tilde{F}_t$.

Then ${\mathcal B}+\nu \tilde{{\mathcal T}}$ generates a unique local
uniformly Lipschitz flow on
$T\mathcal{D}_{\mu}^s$ for $\nu \ge 0$, and the integral curves
$\eta^\nu(t)$ with
$\eta^\nu(0)=e$ extend for a fixed time $\tau >0$ independent of $\nu$ and are
unique.
Further,
$$ \lim_{\nu \rightarrow 0} \eta^\nu(t) = \eta^0(t)$$
for each $t$, $0\le t < \tau$, the limit being in the $H^s$ topology,
$s>(n/2)+1+2\sigma$.
In particular, this holds for $\sigma=2$, and ${\mathcal T}=-\triangle_r$.
\end{prop}

By inverting $(1-\alpha^2\triangle_r)$ in (\ref{avg_euler}) we see that
the dissipation is exactly of the form $-\nu \triangle_r u$, and that this
operator with Dirichlet boundary data generates a strongly continuous 
semi-group.  Setting the nonlinear operator ${\mathcal B}$ to equal  the
geodesic spray of (\ref{avg_euler}), which is $C^\infty$ by Theorem
\ref{Thm2}, we have proven the following:
\begin{cor}\label{Cor3}
For $s>{\frac{n}{2}}+5$ and boundaryless manifolds $M$, solutions to
$$
\begin{array}{c}
\partial_t(1-\alpha^2\triangle_r)u - \nu (1-\alpha^2 \triangle_r)\triangle_r u 
+ \nabla_u(1-\alpha^2\triangle_r)u
-\alpha^2 (\nabla u)^t \cdot \triangle_r u = -\text{\text{\rm grad }}p,\\
\text{\text{\rm div }} u =0, \ \ u(0)=u_0,
\end{array}
$$
converge regularly in $H^s$ to solutions of the inviscid equation with
$\nu=0$.
\end{cor}

Theorem \ref{Thm2} also provides interesting geometric corollaries.
We define the Riemannian exponential map 
${\mathcal E}\text{xp}_e : T_e{\mathcal D}_{\mu,D}^s \rightarrow 
{\mathcal D}_{\mu,D}^s$ of the right invariant metric $\langle \cdot,
\cdot \rangle$ by ${\mathcal E}\text{xp}_e(tu) = \eta(t)$, where 
$t>0$ is sufficiently small, and $\eta(t)$ is the geodesic curve on
${\mathcal D}_{\mu,D}^s$ emanating from $e$ with initial velocity $u$.
Because the above theorem guarantees that geodesics of $\langle \cdot,
\cdot \rangle$ have $C^\infty$ dependence on initial data, 
${\mathcal E}\text{xp}_e$ is well defined, satisfies 
${\mathcal E}\text{xp}_e(0)=e$, and so by the inverse function theorem
we obtain
\begin{cor} \label{Cor4}
For $s>(n/2)+1$, the Riemannian exponential map
${\mathcal E}{\rm xp}_e: T_e{\mathcal D}^s_{\mu,D} \rightarrow 
{\mathcal D}^s_{\mu,D}$
is a local diffeomorphism, and two elements $\eta_1$ and $\eta_2$ of 
${\mathcal D}^s_{\mu,D}$ that are in a sufficiently
small neighborhood of $e$ can be connected by a unique geodesic of $\langle
\cdot, \cdot \rangle$ in ${\mathcal D}^s_{\mu,D}$.
\end{cor}

Note that for the $L^2$ right invariant metric on ${\mathcal D}_\mu^s$
whose geodesic flow gives solutions to the Euler equations, the analogous
local result was obtained by Ebin-Marsden \cite{EM}, but
Shnirelman \cite{Shn} has shown that this local result does not hold globally.
Namely, when $M$ is the unit cube in ${\mathbb R}^3$, he proved the
existence of fluid configurations
which cannot be connected to the identity by an energy minimizing curve.
This has motivated the construction of generalized flows;
Brenier \cite{B} has recently constructed Young 
measure-valued flows that are both Lagrangian and Eulerian in character,
and which give weak solutions to the Euler equations in the sense of
connecting any two fluid configurations (again on the unit cube in ${\mathbb R}
^3$).  The construction of such weak solutions for the weak form of
(\ref{avg_euler}), given on the flat three-torus ${\mathbb T}^3$ by
\begin{align*}
&\int_0^T \int_{{\mathbb T}^3}\left\{ -u \cdot \partial_t \phi -u\otimes u :
\nabla \phi \right.\\
& \qquad \qquad \left.
+\alpha^2\left[\nabla u \cdot \nabla u^t + \nabla u \cdot \nabla u -
\nabla u^t \cdot \nabla u\right] : \nabla(1-\alpha^2\triangle)^{-1}\phi
\right\} dx dt =0
\end{align*}
for all $\phi \in C^\infty([0,T]\times {\mathbb T}^3)$ with $\text{div }\phi
=0$,
is the subject of ongoing research.  In this setting, one generates weak
solutions whose distributional derivatives are Young measures.

The last corollary of Theorem \ref{Thm2} which we shall state concerns
the behavior of the exponential map.
Note that while the group exponential map is only $C^0$ and does not 
cover a neighborhood of the identity, the Riemannian exponential map on
${\mathcal D}_{\mu,D}^s$ is smooth by Theorem \ref{Thm2}, so that in
conjunction with the fact that the right multiplication map is $C^\infty$,
the topological group ${\mathcal D}_{\mu,D}^s$ looks very much like a 
Lie group.
As a consequence of the smoothness of ${\mathcal E}\text{xp}_e$ and
the proof of Theorem 12.1 in \cite{EM}, geodesics
of $\langle\cdot, \cdot \rangle$, which are the solutions of (\ref{avg_euler})
with $\nu=0$, instantly inherit the regularity of the initial data.  
Thus,
\begin{cor}\label{Cor5}
For $s> \frac{n}{2}+1$, let $\eta(t)$ be a geodesic of the right invariant 
metric $\langle \cdot, \cdot \rangle$ on ${\mathcal D}^s_{\mu,D}$, i.e.
$\partial_t\eta(t,x) = u(t,\eta(t,x))$ and $u(t,x)$ is the unique solution
of (\ref{avg_euler}) with $\nu=0$.
If $\eta(0) \in {\mathcal D}_{\mu,D}^{s+k}$  and
${\dot{\eta}}(0) \in T_{\eta(0)}{\mathcal D}_{\mu,D}^{s+k}$
for $0 \le k \le \infty$, then $\eta(t)$ is $H^{s+k}$ for all $t \in I$.
\end{cor}

Our final theorem is geometric and concerns the existence of the weak
Levi-Civita covariant derivative on ${\mathcal D}_{\mu,D}^s$ of the
the weak right invariant metric $\langle \cdot, \cdot \rangle$, as well
as its Riemannian curvature operator.  

Because the metric $\langle \cdot, \cdot\rangle$ is equivalent to an
$H^1$ metric by Korn's inequality, it induces a weak topology relative
to the strong $H^s$ topology, $s>{\frac{n}{2}}+1$, of ${\mathcal D}_{\mu,D}^s$.
In general, there does not exist a weak covariant derivative operator 
associated to a weak metric, nor a bounded Riemannian curvature operator.
Thanks to Theorem \ref{Thm2}, however, these structures do indeed exist.

\begin{thm}\label{Thm3}
Extending $X_\eta, Y_\eta, Z_\eta \in T_\eta {\mathcal D}^s_{\mu,D}$, 
$\eta \in {\mathcal D}_{\mu,D}^s$, to smooth vector
fields $X,Y,Z$ on ${\mathcal D}^s_{\mu,D}$, there exists a right invariant 
unique
Levi-Civita covariant
derivative $\tilde \nabla$ of $\langle \cdot, \cdot\rangle$ on
${\mathcal D}^s_{\mu,D}$
given by
\begin{align*}
\tilde \nabla_XY(\eta)& = \Bigl\{ {\mathcal P}_e \circ  \bigl[
\partial_t(Y_\eta\circ \eta^{-1}) +
       \nabla_{X_\eta\circ\eta^{-1}}(Y_\eta\circ \eta^{-1}) \\
&
+ \frac{1}{2}({\mathfrak U}(X_\eta \circ \eta^{-1}, Y_\eta \circ \eta^{-1})
+{\mathfrak R}(X_\eta\circ \eta^{-1}, Y_\eta\circ \eta^{-1}))
\bigr] \Bigr\} \circ \eta,
\end{align*}
where ${\mathfrak U}$ and ${\mathfrak R}$ are given by polarization of
the operators ${\mathcal U}$ and ${\mathcal R}$, respectively, defined by
\begin{align*}
{\mathcal U}(u) =& \alpha^2(1-\alpha^2{\mathcal L})^{-1}\bigl\{
\text{\rm div}\left[ \nabla u \cdot
\nabla u^t + \nabla u \cdot \nabla u - \nabla u^t \cdot \nabla u\right]
+\text{\rm \text{grad }Tr}(\nabla u \cdot \nabla u)\bigr\} \\
{\mathcal R}(u) =& \alpha^2(1-\alpha^2{\mathcal L})^{-1} \bigr\{  \text{\rm Tr}
\left[ \nabla \left( R(u,\cdot)u \right) +R(u,\cdot) \nabla u +
R(\nabla u, \cdot)u \right] \\
&\qquad\qquad\qquad  + \text{\rm \text{grad }Ric}(u,u)
-(\nabla_u\text{\rm Ric}) \cdot u + \nabla u^t \cdot \text{\rm Ric}(u)
\bigr\}  ,
\end{align*}
and  where for $r\ge 1$, ${\mathcal P}_e: H^r(TM)\cap H^1_0(TM) 
\rightarrow V^r_\mu$  is the
$\langle \cdot, \cdot\rangle_e$-orthogonal projection given by
$${\mathcal P}_e(F) = v$$
where $v\in {\mathcal V}^r_\mu$ is the unique solution of the Stokes problem
$$
\begin{array}{c}
(1-{\mathcal L})v + \text{\rm grad }p=(1-{\mathcal L})F,\\
\text{\rm div }v=0, \ \ v=0 \text{ \rm on } {\partial M}.
\end{array}
$$
For right-invariant vector fields $X, Y$ on ${\mathcal D}^s_{\mu,D}$ which
are completely determined by there value at the identity $X_e, Y_e$,
$$
\tilde \nabla_{X}Y(e) = {\mathcal P}_e  \circ  \Bigl[
       \nabla_{X_e}Y_e + \frac{1}{2}\bigl({\mathfrak U}(X_e, Y_e)
+ {\mathfrak R}(X_e,Y_e)\bigr)
\Bigr] .
$$

Finally,  define
the weak Riemannian curvature tensor
$$ \tilde R_\eta : \left[T_\eta {\mathcal D}^s_{\mu,D} \right]^3 \rightarrow
T_\eta {\mathcal D}^s_{\mu,D}$$
by
$$ \tilde R_\eta(X_\eta,Y_\eta)Z_\eta= \left( \tilde \nabla_Y \tilde \nabla_X Z
\right)_\eta - \left( \tilde \nabla_X \tilde \nabla_Y Z \right)_\eta
+ \left( \tilde \nabla_{[X,Y]} Z \right)_\eta, \ \ \eta\in
{\mathcal D}_{\mu,D}^s.$$
Then for $s>(n/2)+2$, $\tilde R$ is right invariant and continuous
in the $H^s$ topology.
\end{thm}

Since the weak curvature operator $\tilde R$ is bounded in $H^s$ for
$s>{\frac{n}{2}}+2$, the fundamental existence and
uniqueness theorem for ordinary differential equations provides us with the 
following:
\begin{cor}\label{Cor6}
For $s>(n/2)+2$ and $y, \dot y \in T_e {\mathcal D}^s_{\mu,D}$, there exists
a unique $H^s$ vector field $Y(t)$ along a geodesic curve $\eta$ of $\langle \cdot,
\cdot \rangle$ which is solution to the Jacobi equation
$$ \tilde\nabla_{\dot\eta} \tilde\nabla_{\dot\eta} Y +
\tilde{R}_\eta (\dot\eta,Y) \dot\eta =0,
\ \ Y(0)=y,\ \tilde\nabla_{\dot\eta}Y(0)=\dot{y} .
$$
\end{cor}

Because the geodesic flow $\eta$ of the right invariant metric $\langle \cdot,
\cdot \rangle$ on ${\mathcal D}_{\mu,D}^s$ is the solution of (\ref{avg_euler})
with $\nu=0$, and since Jacobi's equation is the linearization of the
geodesic flow, Corollary \ref{Cor6} proves existence and uniqueness of
(\ref{avg_euler}), linearized about a solution $u= \partial_t \eta \circ
\eta^{-1}$.  We are thus able to follow Arnold \cite{A}, and study the
Lagrangian stability of our solutions, by studying the curvature of our
infinite-dimensional group.  Positive curvature indicates stable motion,
while negative curvature implies exponential divergence of trajectories.

Since this system, thought of as the averaged Euler equation, averages over
the small-scale fluid motion, one might hope that solutions of (\ref{avg_euler})
might have nicer stability properties than solutions to the Euler equations.
Geometrically, this implies that as $\alpha$ is increased away from zero, the
sectional curvatures which are negative for Euler flow flip sign and become
positive.  Indeed, this seems to be the case; we give a simple example.

We consider periodic two-dimensional motion, so the configuration space
is the group of volume-preserving diffeomorphisms of the two-torus 
${\mathbb T}^2$.
Consider the parallel sinusoidal steady flow
given by the stream function $\xi = \cos (k, x)$
and let $\psi$ be any other vector of the tangent space at $e$, i.e.,
$\psi = \sum x_l e_l$, where $x_{-l} = \bar{x}_l$.
Theorem $3.4$ of \cite{AK} states that the curvature of the group
${\mathcal D}_\mu (T^2)$ in any two-dimensional plane containing
the direction $\xi$ is \emph{non-positive} and is given by
$$
  K_{\xi \psi} = \frac{S}{4} \sum_l a_{k l}^2 |x_l + x_{l+2k} |^2,
$$
where
$a_{k l} = \dfrac{(k \times l)^2}{| k + l |}$,
$k \times l = k_1 l_2 - k_2 l_1$ is the (oriented) area of the
parallelogram spanned by
$k$ and $l$, and $S$ is the area of the torus.
Then, a corollary of this theorem states that the curvature
in the plane defined by the stream functions $\xi = \cos (k,x)$
and $\psi = \cos (l,x)$ is
$$
  K_{\xi \psi} = - (k^2 + l^2) \sin^2 \beta \sin^2 \gamma / 4S,
$$
where $\beta$ is the angle between $k$ and $l$, and $\gamma$ is
the angle between $k+l$ and $k-l$.  Recall that these are the curvatures
with respect to the right invariant $L^2$ metric.

Now using the right invariant metric $\langle \cdot , \cdot \rangle$
on ${\mathcal D}_\mu^s({\mathbb T}^2)$, one can prove the following result:
\footnote{This result was obtained together with Sergey Pekarsky.}
Let $\tilde K(\xi,\psi)$ denote the sectional curvature on
${\mathcal D}_\mu^s({\mathbb T}^2)$ with the right invariant metric
$\langle \cdot, \cdot \rangle$, where $\xi=\cos(k,x)$ and $\psi=\cos(l,x)$.
For $|\epsilon|$ sufficiently small, let $l=k+\epsilon$.  Then
for any $k$, there exists $0 < \alpha_0(k) < 1$, such that for all
$\alpha > \alpha_o(k)$, $\tilde K(\xi, \psi) >0$.

\section{Review of the Hilbert manifold of maps and diffeomorphism groups}
\label{sec1}

Let us briefly recall some facts concerning the geometry of the manifold
of maps between two Riemannian manifolds.  We refer the reader to 
\cite{Pal}, \cite{Eells}, and \cite{Elias} for a comprehensive treatment of
this subject.  Let $(M,g)$ be a $C^\infty$ compact oriented $n$-dimensional
Riemannian manifold  with boundary, and let $(N,h)$ denote a 
$p$-dimensional compact oriented boundaryless Riemannian manifold.
By Sobolev's embedding theorem,
when $s > n/2 + k$, the set of Sobolev mappings $H^s(M,N)$
is a subset of $C^k(M,N)$ with continuous inclusion, and so
for $s>n/2$, an $H^s$-map of $M$ into $N$ is pointwise
well-defined.  Mappings in the space $H^s(M,N)$ are those whose
first $s$ distributional derivatives are square integrable in {\it
any} system of charts covering the two manifolds.

For $s>n/2$, the space $H^s(M,N)$ is a $C^\infty$
differentiable Hilbert manifold.  Let exp$:TN
\rightarrow N$ be the exponential mapping associated with $h$.
Then for each $\phi \in H^s(M,N)$, the
map $\omega_{\operatorname{exp}}:T_\phi H^s(M,N) \rightarrow
H^s(M, N)$ is used to provide a differentiable structure which
is independent of the chosen metric, where
$\omega_{\operatorname{exp}}(v) = \operatorname{exp} \circ v$, and
$T_\phi H^s(M,N) =\{ u \in H^s(M,TN)\ | \ \bar\pi \circ u = \phi\}$, where
$\bar\pi:TN \rightarrow N$.

When ${\partial M} \neq \emptyset$, the set $H^s(M,M)$ is not a smooth
manifold. We can, however, embed $\overline M$ into its double $\tilde M$,
a compact boundaryless manifold of the same dimension, extending the metric
$g$ to $\tilde M$.  Using the above construction, we form the $C^\infty$ manifold
$H^s(M,\tilde{M})$.  Then for $s>(n/2)+1$, the set
\begin{align*}
{\mathcal D}^s = \{ \eta \in H^s(M,\tilde{M}) \ | \ \eta &\text{ is bijective },
\eta^{-1}\in H^s(M,\tilde{M}), \\
& \eta \text{ leaves } {\partial M} \text{ invariant}\}
\end{align*}
is an open subset of $H^s(M,\tilde{M})$.  By choosing a metric on $\tilde M$
for which ${\partial M}$ is a totally geodesic submanifold, the above construction
provides ${\mathcal D}^s$ with a $C^\infty$ differentiable structure (see
\cite{EM} for details).  For each $\eta \in {\mathcal D}^s$, the tangent
space at $\eta$ is given by
\begin{align*}
T_\eta{\mathcal D}^s=\{ u\in H^s(M,TM) \ | \ \pi \circ  u = \eta, \ \ g(u\circ
\eta^{-1},n)=0
\text{ on } {\partial M} \}
\end{align*}
and the vector space $T_e {\mathcal D}^s$ consists of the $H^s$ class vector
fields on $M$ which are tangent to ${\partial M}$.

Let $\mu$ denote the Riemannian volume form on $M$, and let 
$$ \mathcal{D}_\mu^s := \{\eta \in \mathcal{D}^s \mid \eta^*(\mu)=\mu \}$$
be the subset of ${\mathcal D}^s$ whose elements preserve $\mu$.  As
proven in \cite{EM}, the set ${\mathcal D}_\mu^s$ is a $C^\infty$ subgroup of 
$\mathcal{D}^s$ for $s>(n/2)+1$.  We call ${\mathcal D}_\mu^s$ the group
of volume preserving diffeomorphisms of class $H^s$.  The tangent space at
$\eta\in {\mathcal D}_\mu^s$ is given by
\begin{align*}
T_\eta{\mathcal D}_\mu^s=\{ u\in &H^s(M,TM) \ |  \ \pi \circ  u = \eta, \ \ 
g(u\circ \eta^{-1},n)=0
\text{ on } {\partial M},\\
&\text{div}(u \circ \eta^{-1}) =0 \},
\end{align*}
so that the vector space $T_e {\mathcal D}_\mu^s$ consists of divergence-free
$H^s$ class vector fields on $M$ that are tangent to ${\partial M}$.

We have the following standard composition lemma:

\begin{lemma}[$\omega$ and $\alpha$ lemmas]\label{o-lemma}
For $\eta \in {\mathcal D}^s$, right multiplication
$$ R_\eta : {\mathcal D}^s \rightarrow {\mathcal D}^s
\ \ (H^s \rightarrow H^s), \ \zeta \mapsto \zeta \circ \eta, 
\text{ is } C^\infty,$$
and for $\eta\in {\mathcal D}^{s+r}$, left multiplication
$$ L_\eta : {\mathcal D}^s \rightarrow {\mathcal D}^{s}
\ \ (H^s \rightarrow H^{s}), \ \zeta \mapsto \eta \circ \zeta, 
\text{ is } C^r.$$
\end{lemma}

Finally, the inverse map $(\eta \mapsto \eta^{-1}): {\mathcal D}^s 
\rightarrow {\mathcal D}^s$ is only $C^0$  and not even locally Lipschitz
continuous.  Thus, ${\mathcal D}^s$ and ${\mathcal D}_\mu^s$ are not Lie groups,
but are $C^\infty$ topological groups with $C^\infty$ right translation.

\section{Proof of Theorem \ref{Thm1}}
\label{sec2}
\subsection{The Neumann group ${\mathcal D}_{\mu,N}^s$}
We begin by first establishing the result for 
${\mathcal D}_{\mu,N}^s$.  We split the proof into three steps.

\medskip
\noindent
{\bf Step 1.} {\bf Bundles over ${\mathcal D}_\mu^s$ and the
transversal mapping theorem.}

\medskip
\noindent
Recall that a smooth map between Hilbert manifolds $f: {\mathcal M}_1
\rightarrow {\mathcal M}_2$ is transversal to a submanifold ${\mathcal M}_3$
of ${\mathcal M}_2$ if for all $m \in f^{-1}({\mathcal M}_3)$,
$Tf(m)\left(T_m{\mathcal M}_1\right) + T_{f(m)}{\mathcal M}_3
= T_{f(m)}{\mathcal M}_2$.  The transversal mapping theorem asserts that
$f^{-1}({\mathcal M}_3)$ is a submanifold of ${\mathcal M}_1$ if $f$ is
transversal to ${\mathcal M}_3$.

Let us define the following infinite dimensional vector bundles over
${\mathcal D}_\mu^s$:

\begin{align*}
{\mathcal F} & = H^{s-{\frac{3}{2}}}_\eta(TM|\partial M) 
             \downarrow {\mathcal D}_\mu^s, \\
{\mathcal E} & = H^{s-{\frac{3}{2}}}_\eta (T\partial M) 
             \downarrow {\mathcal D}_\mu^s, \\
{\mathcal G} & = \left[H^{s-{\frac{3}{2}}}_\eta (TM|\partial M)^*
                 \otimes H^{s-{\frac{3}{2}}}_\eta (T\partial M)\right] 
             \downarrow {\mathcal D}_\mu^s.
\end{align*}

For $x \in {\partial M}$, let $\Pi_x: T_xM \rightarrow T_x{\partial M}$
be the $g$-orthogonal projector, and define the section $\Pi:
{\mathcal D}_\mu^s \rightarrow {\mathcal G}$ pointwise by
$\Pi(\eta)(x) = \Pi_{\eta(x)}$, so that for all $\eta \in
{\mathcal D}_\mu^s$, $\Pi(\eta): H^{s-3/2}_\eta(TM|{\partial M})
\rightarrow H^{s-3/2}_\eta(T{\partial M})$.  For $n \in H^{s-1/2}(N)$,
define the section of ${\mathcal F}$, $h_n: {\mathcal D}_\mu^s \rightarrow
{\mathcal F}$, by
$$h_n(\eta) = T\eta|_{\partial M} \cdot n.$$
Finally, let $f_n: {\mathcal D}_\mu^s \rightarrow {\mathcal E}$ denote the
section of ${\mathcal E}$ which is given by
$$ f_n= \Pi \circ h_n.$$
Then, the set ${\mathcal D}^s_{\mu,N}$ is the inverse image of $f_n$ acting on
the zero section of ${\mathcal E}$.

\begin{lemma}\label{smooth2}
The map $f_n: {\mathcal D}^s \rightarrow {\mathcal E}$ is $C^\infty$.
\end{lemma}
\begin{proof}
This follows from Lemma \ref{smooth1}, the trace theorem, and the fact that
$\Pi$ is smooth, as $g$ and ${\partial M}$ are $C^\infty$.  \end{proof}

Hence, by the transversal mapping theorem, to show that ${\mathcal D}^s_{\mu,N}$
is a $C^\infty$ subgroup of ${\mathcal D}_\mu^s$, we shall prove that $f_n$ is a
surjection;  this will provide ${\mathcal D}_{\mu,N}^s$ with smooth 
differentiable structure.  That ${\mathcal D}^s_{\mu,N}$ is a $C^\infty$
subgroup then follows from the fact that
${\mathcal D}^s_{\mu,N}$ is trivially closed under right
composition.

\medskip
\noindent
{\bf Step 2.} {\bf The covariant derivative of $f_n$.}

\medskip
\noindent
We use the symbol $\overline \nabla$ to denote the weak Levi-Civita
covariant derivative on sections of ${\mathcal F}$ and ${\mathcal G}$ (as 
obtained in Lemma \ref{smooth1}).  Following the methodology of 
Lemma \ref{smooth1},
we compute that for all $\eta \in {\mathcal D}_\mu^s$ and $u \in
T_\eta{\mathcal D}_\mu^s$, $\overline \nabla_u h_n(\eta) \in {\mathcal F}
_\eta = H^{s-3/2}_\eta (TM|{\partial M})$ is given by
$$ \overline \nabla _u h_n (\eta) = \nabla_n u,$$
where $\nabla$ denotes the Levi-Civita covariant derivative in 
$\eta^*(TM)$.

Next, we compute the covariant derivative of the section $\Pi$ of 
${\mathcal G}$.  We shall denote the metric tensor $g$ evaluated at the point
$\eta(x)$ by $g_{\eta(x)}$.  Using the fact that $g$ is covariantly constant,
and letting $( \cdot  )^{\rm tan}$ denote the tangential component of a
mapping $v:{\partial M} \rightarrow TM|{\partial M}$, we have that

\begin{eqnarray}
&&g_{\eta(x)} \left([\nabla _u \Pi_{\eta(x)}] \cdot v(x) , z(x)\right)
= - g_{\eta(x)}\left( (\nabla_u v(x))^{{\rm tan}}, z(x) \right)  \nonumber\\
&& \qquad \qquad  - g_{\eta(x)} \left( (\nabla_u z(x))^{{\rm tan}},
v(x) \right) - u \left[ g_{\eta(x)}(v^{{\rm tan}} (x), z^{{\rm tan}}(x))
\right] \label{1}
\end{eqnarray}
where we use the notation: $u[f] = df \cdot u$ for any function 
$f \in C^1(M)$.  It is clear
that the operator $\nabla_u \Pi_\eta$ is self-adjoint with respect to $g$.
By definition of the $g$-orthogonal projector $\Pi_{\eta(x)}$, we see that
for all $x \in {\partial M}$,
$$ g_{\eta(x)}\left( \Pi_{\eta(x)} \cdot w(x), \nu(x)\right) =0, \ \
\forall \ w\in {\mathcal F}_\eta, \nu \in H^{s-3/2}_\eta(N),$$
so that setting the map $v$ in equation (\ref{1}) equal to the mapping 
$\nu$, and noting that the covariant derivative $\overline \nabla$ on 
${\mathcal G}$ is the functorial lift of $\nabla$, we obtain the formula
$$ \left[ \overline \nabla_u \Pi(\eta) \right] (\nu) = - (\nabla_u \nu)
^{{\rm tan}} = S_\nu(u).$$
It follows that for all $\eta\in f_n^{-1}(0)$,
\begin{eqnarray*}
\overline\nabla_u f_n(\eta) &=& 
\overline\nabla_u \Pi_\eta \cdot h(\eta) + \Pi_\eta
\overline\nabla_u h(\eta) \\
&=& S_\nu (u) + (\nabla_n u)^{{\rm tan}} \in {\mathcal E}_\eta,
\end{eqnarray*}
where $\nu = T\eta|_{\partial M} \cdot n \in H^{s-\frac{3}{2}}_\eta(N)$.

\medskip 
\noindent
{\bf Step 3.} {\bf $f_n$ is a surjection.}

\medskip 
\noindent
It remains to show that for all $\eta \in f_n^{-1}(0)$, $\overline\nabla f_n(\eta):
T_\eta {\mathcal D}_\mu^s \rightarrow {\mathcal E}_\eta$ is onto.  Because
right translation on ${\mathcal D}_\mu^s$ is a smooth operation, it suffices
to find $u\in T_e {\mathcal D}_\mu^s$ such that $\overline \nabla_u f_n(e)=w$
for any $w \in H^{s-3/2}(T{\partial M})$.  To do so, we shall solve the
following elliptic boundary value problem: Find $(u,p) \in 
T_e{\mathcal D}_\mu^s \times H^{s-1}(M)/{\mathbb R}$ such that
\begin{equation}\label{bvp1}
\begin{array}{c}
(1-\triangle_r) u + \text{grad} \ p = F, \ \ \text{div }u=0 \text{ in } M, \\
g(u,n)=0, \ \ (\nabla_n u)^{\rm tan} + S_n(u) = w \text{ on } {\partial M},
\end{array}
\end{equation}
where $F \in H^{s-2}(TM), w \in H^{s-3/2}(T{\partial M}), n \in H^{s-1/2}(N)$.

We first define the space
$$H^1_A(TM) = \{v \in H^1(TM) | \text{div}\ v=0 \text{ and } g(u,n)=0
\text{ on } {\partial M}\},$$
and establish the existence of a unique weak solution $u\in H^1_A(TM)$ to 
(\ref{bvp1}).  Let $B: H^1_A(TM) \times H^1_A(TM) \rightarrow {\mathbb R}$
be the bilinear form given by
$$B(u,v) = \int_M \left[  g(u,v) + 
2\bar g(\text{Def } u,\text{Def }v)\right] \mu .$$
$B$ is symmetric and by Korn's inequality, which states that
$|u|_1 \le C |\text{Def }u|_0 + C |u|_0$ 
(see, for example, \cite{T} Corollary 12.3),  there exists $\beta >0$ such that
$\beta | u|_1 \le B(u,u)$; hence, $B$ is coercive with respect to $H^1_A(TM)$.
Let ${\mathfrak F} : H^1_A(TM) \rightarrow {\mathbb R}$ be given by
${\mathfrak F}(v) = \int_M g(F,v)\mu + \int_{\partial M} g(w,v) \mu_\partial$.
By the trace theorem, $|\int_{\partial M} g(w,v) \mu_\partial | \le C 
|w|_{L^2(T{\partial M})} \ |v|_1$, so that together with the Cauchy-Schwartz
inequality and the embedding $H^1 \hookrightarrow L^2$, we see that 
${\mathfrak F} \in H^1_A(TM)^*$.  Hence, by the
Lax-Milgram theorem, their exists a unique $u \in H^1_A(TM)$ satisfying
$B(u,v) = {\mathfrak F}(v)$ for all $v \in H^1_A(TM)$.  This, in turn,
uniquely determines $p\in L^2(M)/{\mathbb R}$, as the
solution of $B(u,v) - {\mathfrak F}(v)= \int_M p \cdot \text{div }v \mu$ for all
$v \in H^1(TM)$ that satisfy $g(u,n)=0$ on ${\partial M}$. We have thus obtained
a unique weak solution $(u,p) \in H^1_A(TM) \times L^2(M)/{\mathbb R}$ of the
boundary value problem (\ref{bvp1}).

Now, since 
$$2\text{Def}^*\text{Def }u = -2\text{Div}\text{Def }u 
  =-\triangle u - 2\text{Ric}(u),$$
we see that if $u\in H^2(TM)\cap H^1_A(TM)$ satisfies
$$B(u,v) = {\mathfrak F}(v), \ \ \forall \ v\in H^1_A(TM),$$
then u is a solution of (\ref{bvp1}).  We shall use an elliptic regularity
argument to prove that $u$ is in fact a classical $H^s$ solution of 
(\ref{bvp1}).

Let $(U,\phi)$ coordinate chart on $\overline M$, and $\chi \in C^\infty_0 (U)$.
Since $(1-\triangle_r) (\chi u) = \chi( (1-\triangle_r)u) + 
[(1-\triangle_r), \chi] u$, and since $[(1-\triangle_r), \chi] u$ is a 
first-order
differential operator, our elliptic regularization of $u$ can be localized
to the chart $U$.  We can assume that $U$ intersects ${\partial M}$, for 
otherwise, standard interior regularity estimates can be applied.  Let
$x^i$ denote the coordinates on $U$ and set $\partial_i =\partial / \partial x^i$.
We may express the Hodge Laplacian $\triangle$ on $U$ as
$$ \triangle u = \triangle_{\rm loc} + Y(u),$$
where $\triangle_{\rm loc}= g^{ij}(x) \partial_i \partial_j u$, and $Y$ is a 
first order differential operator.

We consider the boundary value problem in $U$ given by
\begin{equation} \nonumber
\begin{array}{c}
(1-\triangle_{\rm loc}) u + \text{grad }p=F, \ \ \text{div }u=\rho \text{ in } U,\\
B_0(u)=0, {B}_1(u)=w \text{ on } {\partial U},
\end{array}
\end{equation}
where $B_0(u)=g(u,n)$, $B_1(u) = 2[ (Du + Du^t) \cdot n]^{\rm tan}$, and
$Du\cdot n = \partial_j u^i g^j_k n^k$.
Applying induction  to the usual difference quotient argument
(see, for example, \cite{T}) yields the elliptic estimate
$$|u|_s + |p|_{s-1} \le C\left( |F|_{s-2} +|\rho|_{s-1} +|B_0(u)|_{s-1/2}+
|{B}_1(u)|_{s-3/2}\right).$$
Hence,  the operator
$\overline {\mathcal L}: H^s \cap H^1_A(TU) \rightarrow H^{s-2}(TU)
\oplus H^{s-1}(U)$$\oplus$$H^{s-1/2}(T{\partial U})$
$\oplus H^{s-3/2}(T{\partial U})$ given by
$$ \overline {\mathcal L} u = \left( (1-\triangle_{\rm loc}u),\text{div }u,B_0(u),
{B}_1(u) \right)= (F,\rho,0,w)$$
has closed range, and since its adjoint has a trivial kernel, 
$\overline {\mathcal L}$  is an isomorphism
(see also \cite{MPS} for an alternative proof that $\overline{\mathcal L}$ 
is an isomorphism).

A simple computation verifies that along ${\partial M}$, 
$$ 2[\text{Def }u \cdot n]^{\rm tan} =
(\nabla_nu)^{\rm tan}+S_n(u) \ \ \forall \ u\in H^1_A(TM),$$
so that on
${\partial U}$, $[\text{Def }u\cdot n]^{\rm tan}$ differs from $B_1(u)$ by a linear
combination of $C^\infty$ Christoffel maps, and we shall denote this difference by
$\Gamma(u)$.  Hence, the operator 
$\mathcal L: H^s \cap H^1_A(TU) \rightarrow H^{s-2}(TU)
\oplus H^{s-1}(U) \oplus H^{s-1/2}(T{\partial U})
\oplus H^{s-3/2}(T{\partial U})$ given by
$${\mathcal L}u =\left((1-\triangle_r)u,\text{div }u, B_0(u),
(\nabla_nu)^{\rm tan}+S_n(u)\right)$$
differs from $\overline{\mathcal L} u$ by the operator ${\mathcal K}u
=(Y(u)+\text{Ric}(u),0,0,\Gamma(u))$ which is compact by Rellich's theorem.
Therefore, ${\mathcal L}$ has index $0$ and trivial kernel, and is thus an
isomorphism, which concludes that ${\mathcal D}_{\mu,N}^s$ is a 
$C^\infty$ subgroup of ${\mathcal D}_\mu^s$.

With an almost trivial modification, ${\mathcal D}_{N}^s$ is a $C^\infty$
subgroup of ${\mathcal D}^s$.  To see this,
we redefine the vector bundles ${\mathcal E}, {\mathcal F}, {\mathcal G}$ to
have ${\mathcal D}^s$ as base manifold rather than ${\mathcal D}_\mu^s$,
and we redefine the space $H^1_A(TM)$, removing the divergence-free 
constraint.  In this case, 
$$2 \text{Def}^* \text{Def } u = -2\text{Div}\text{Def }u = -(\triangle + 2\text{Ric}
+\text{grad }\text{div})u,$$
so to establish that $f_n$ is a surjection, we solve the following boundary value problem:
For $F \in H^{s-2}(TM)$, $w\in H^{s-3/2}(T{\partial M})$ and 
$n \in H^{s-1/2}(N)$, find $u\in T_e{\mathcal D}^s$ satisfying
\begin{equation}\label{bvp1a}
\begin{array}{c}
[1-(\triangle_r +\text{grad }\text{div})]u = F \text{ in } M\\
g(u,n)=0, \ \ (\nabla_nu)^{\rm tan} + S_n(u)=w \text{ on } {\partial M}.
\end{array}
\end{equation}

A weak solution in $H^1_A(TM)$ is obtained using the Lax-Milgram theorem just
as in Step 3 above.
Up to a compact operator, this is precisely the elliptic system studied
in (\cite{F}), wherein existence and uniqueness of classical $H^s$ solutions
is established.  Since 
modification of an elliptic operator by lower-order terms does not change
its index, we have existence of $u \in T_e {\mathcal D}_\mu^s$ solving 
(\ref{bvp1a}), and this completes the argument for the subgroup
${\mathcal D}_{N}^s$.

\subsection{The mixed group ${\mathcal D}_{\mu,mix}^s$}

We shall follow the three step proof above, keeping the same notation.
 
\medskip
\noindent
{\bf Step 1. Bundles over ${\mathcal D}_\mu^s$ and the inverse function
theorem.}

\medskip
\noindent
We modify the vector bundles ${\mathcal F}$, ${\mathcal E}$, and ${\mathcal G}$
as follows:
\begin{align*}
{\mathcal F} & = H^{s-{\frac{3}{2}}}_\eta(TM|\Gamma_2) 
             \downarrow {\mathcal D}_\mu^s, \\
{\mathcal E} & = H^{s-{\frac{3}{2}}}_\eta (T\Gamma_2) 
             \downarrow {\mathcal D}_\mu^s, \\
{\mathcal G} & = \left[H^{s-{\frac{3}{2}}}_\eta (TM|\Gamma_2)^*
                 \otimes H^{s-{\frac{3}{2}}}_\eta (T\Gamma_2)\right] 
             \downarrow {\mathcal D}_\mu^s.
\end{align*}
For $n \in H^{s-1/2}(N|{\Gamma_2})$, define
$\bar f_n:{\mathcal D}_\mu^s \rightarrow {\mathcal D}^{s-1/2}(\Gamma_1)
\times {\mathcal E}$ by
$$\bar f_n (\eta) = \left[ \eta|_{\Gamma_1}, f_n(\eta)\right]
=\left[ \eta|_{\Gamma_1}, \Pi(\eta)\circ (T\eta|_{\Gamma_2}\cdot n)\right].$$
The trace theorem together with Lemma \ref{smooth2} ensures that $\bar f_n$ is
$C^\infty$.  Since ${\mathcal D}^S_{\mu,mix}$ $=$ $\bar f_n^{-1}(e,0)$, we must prove
that $\bar f_n$ is a surjection, in order to show that ${\mathcal D}^S_{\mu,mix}$
is a submanifold of ${\mathcal D}_\mu^s$.  Again, it is clear that the
set ${\mathcal D}^S_{\mu,mix}$ is closed under right composition.

\medskip 
\noindent
{\bf Step 2. Computing the tangent map of $\bar f_n$.}

\medskip 
\noindent
Step 2 of the case ${\mathcal D}_{\mu,N}^s$ shows that for any $u \in 
T_\eta {\mathcal D}_\mu^s$,
$$ \overline \nabla_u f_n = S_\nu(u) + (\nabla_nu)^{\rm tan} \in {\mathcal E}
_\eta, \ \ \nu = T\eta|_{\Gamma_2} \cdot n \in H^{s-3/2}_\eta(N|\Gamma_2).$$

Now $\overline \nabla_u f_n$ is the vertical component of $Tf_n \cdot u$,
the $T{\mathcal E}$-valued image of $u$ under the tangent mapping $Tf_n$.
Letting ${\mathcal H}\subset T{\mathcal E}$ denote the connection associated
with the Levi-Civita covariant derivative $\overline \nabla$ (see Step 1 above),
we have the
local decomposition $Tf_n \cdot u = \overline \nabla_u f_n - \omega_{\mathcal H}
(u) \cdot f_n$, where $\omega_{\mathcal H}$ is the local connection $1$-form
on ${\mathcal E}$ associated with the horizontal distribution ${\mathcal H}$.
Then, 
$$T\bar f_n(\eta) \cdot u = \left(u|_{\Gamma_1}, 
\overline \nabla _u f_n(\eta) - \omega_{\mathcal H}(u)\cdot f_n (\eta)\right.),
 \ \ u \in T_\eta{\mathcal D}_\mu^s.$$

\medskip 
\noindent
{\bf Step 3. $\bar f_n$ is a surjection.}  

\medskip 
\noindent
It suffices to prove that for all
$(\psi, w) \in H^{s-1/2}(TM|\Gamma_1) \times {\mathcal E}_e$, there exists $u\in T_e
{\mathcal D}_\mu^s$ such that
\begin{align*}
u & = \psi \text{ on } \Gamma_1 \\
(\nabla_nu)^{\rm tan} + S_n(u) & = w \text{ on } \Gamma_2,
\end{align*}
and to do so, we shall follow
Step 3 for the case of ${\mathcal D}_{\mu,N}^s$, and obtain $u$ as the 
solution of
\begin{equation}\label{bvp2}
\begin{array}{c}
(1-\triangle_r) u + \text{grad} \ p = F, \ \ \text{div }u=0, \text{ in } M, \\
u=\psi \text{ on } \Gamma_1,\\
g(u,n)=0, \ \ (\nabla_n u)^{\rm tan} + S_n(u) = w \text{ on } {\Gamma_2}.
\end{array}
\end{equation}
It suffices to consider the homogeneous boundary condition $u=0$ on $\Gamma_1$.

To obtain a weak solution to (\ref{bvp2}), we define
$$H^1_A(TM) = \{v \in H^1(TM) | \text{div}\ v=0, \  g(u,n)=0
\text{ on } \Gamma_2 \text{ and } u=0 \text{ on } \Gamma_1\},$$
and again consider the bilinear form $B: H^1_A(TM) \times H^1_A(TM) \rightarrow
{\mathbb R}$ given by
$$B(u,v)=\int_M\left[g(u,v)+2\bar g(\text{Def }u,\text{Def }v)\right]\mu .$$
We define ${\mathfrak F}: H^1_A(TM) \rightarrow {\mathbb R}$ by
${\mathfrak F}(v) = \int_M g(F,v)\mu + \int_{\Gamma_2} g(w,v) \mu_\partial$.
The argument we gave in Step 3 of the case ${\mathcal D}_{\mu,N}^s$ shows that
there exists a unique solution $u\in H^1_A(TM)$ satisfying $B(u,v)={\mathfrak F}(v)$
for all $v\in H^1_A(TM)$.

Now, if $u\in H^2(TM) \cap H^1_A(TM)$ satisfies
$B(u,v) = {\mathfrak F}(v)$ for all $v \in H^1_A(TM)$, then $u$ is a solution of
the mixed problem (\ref{bvp2}) for which elliptic regularity is slightly more
subtle than for the Neumann problem.  In particular, the identical argument 
which we used for {\it that} problem provides the $H^s$ class regularity of $u$
on $M/(\Gamma_1 \cap \Gamma_2)$; after all, the boundary 
conditions on both $\Gamma_1$ and $\Gamma_2$ are elliptic in the sense of
Agmon-Douglis-Nirenberg as the Complementary Condition is satisfied (see
\cite{ADN}, and see \cite{SS} for an alternative method).  The fact that
${\partial M}$ is $C^\infty$ and  that
${\partial M}= \Gamma_1 \cup \Gamma_2$ gives the regularity of the solution 
on $\overline M$ (see, for example, Fichera \cite{F}, pages 377 and 385).
Hence, our argument in Step 3 for the subgroup ${\mathcal D}_{\mu,N}^s$ given
above yields a unique solution $u\in H^s(TM)\cap H^1_A(TM)$ of
(\ref{bvp2}), and thus concludes the proof that ${\mathcal D}_{\mu,mix}^s$
is a $C^\infty$ subgroup of ${\mathcal D}_\mu^s$.

Just as we proved that ${\mathcal D}_{N}^s$ is a subgroup of ${\mathcal D}^s$
by a minor modification of the argument for the case ${\mathcal D}_{\mu,N}^s$,
we easily obtain that ${\mathcal D}_{mix}^s$ is also a $C^\infty$ subgroup 
of ${\mathcal D}^s$.

\subsection{The subgroup ${\mathcal D}_{\mu,D}^s$}
This case was studied by Ebin-Marsden \cite{EM} using a different approach.
By setting $\Gamma_2=\emptyset$ above, we immediately prove that
${\mathcal D}_{\mu,D}^s$ is a $C^\infty$ subgroup of ${\mathcal D}_\mu^s$ and
that ${\mathcal D}_{D}^s$ is a $C^\infty$ subgroup of ${\mathcal D}^s$,
with the appropriate tangent spaces at the identity.

\medskip
\noindent
This concludes the proof of Theorem \ref{Thm1}.
\subsection{The group exponential map.}
Let  ${\mathfrak G}^s$ denote either of the groups ${\mathcal D}_{D}^s$, 
${\mathcal D}_{N}^s$, or ${\mathcal D}_{mix}^s$, and similarly, let
${\mathfrak G}^s_\mu$ denote either of the groups
${\mathcal D}_{\mu,D}^s$, ${\mathcal D}_{\mu,N}^s$, or
${\mathcal D}_{\mu,mix}^s$. 

\begin{cor}\label{Exp}
Let $V \in T_e {\mathfrak G}^s$, and let $\eta_t$
be its flow, $(d/dt)\eta_t = V \circ \eta_t$.  Then, for $s>(n/2)+2$, 
$\eta_t$ is a one parameter subgroup of ${\mathcal G}^s$, and the group exponential
map Exp$:T_e{\mathfrak G}^s \rightarrow {\mathfrak G}^s$ given by $V \mapsto \eta_1$
is continuous but not continuously differentiable, while the curve $t \mapsto \eta_t$
is $C^1$.  This holds for ${\mathfrak G}^s_\mu$ as well.
\end{cor}
\begin{proof}
The result follows from (\cite{EM}, Theorems 3.1 and 6.3).
\end{proof}

\subsection{Further remarks on diffeomorphism subgroups}

The existence of the above $C^\infty$ subgroups follows from the existence,
uniqueness, and regularity of solutions to certain elliptic boundary value 
problems.  

\medskip
\noindent
This methodology allows to prove directly that
for $s>(n/2)+1$, ${\mathcal D}_{\mu,mix}^s$ is a $C^\infty$ subgroup of 
${\mathcal D}^s$.

\medskip
\noindent
We need only modify the map $\bar f_n$ given in Step 3 above as follows:
For $n \in H^{s-1/2}(N|\Gamma_2)$ and $\mu$ the Riemannian
volume form on $M$, define
$\bar f_{n,\mu}:{\mathcal D}_\mu^s \rightarrow \Lambda^3(M) \times 
{\mathcal D}^{s-1/2}(\Gamma_1) \times {\mathcal E}$ by
$$\bar f_{n,\mu} (\eta) 
=\left[\eta^*(\mu), \eta|_{\Gamma_1}, 
\Pi(\eta)\circ (T\eta|_{\Gamma_2}\cdot n)\right].$$
Again $\bar f_{n,\mu}$ is $C^\infty$, and following the notation of Step 2,
we easily compute that
$$T\bar f_{n,\mu}(\eta) \cdot u = \left(\text{div }(u \circ \eta^{-1}),
u|_{\Gamma_1}, 
\overline \nabla _u f_n(\eta) - \omega_{\mathcal H}(u)\cdot f_n (\eta)\right.), \ \ u \in T_\eta{\mathcal D}_\mu^s.$$

Finally, the modification to Step 3 consists
of obtaining a solution $u\in T_e{\mathcal D}_\mu^s$ satisfying the boundary
value problem
\begin{equation}\nonumber
\begin{array}{c}
(1-\triangle_r) u + \text{grad} \ p = F, \ \ \text{div }u=q, \text{ in } M, \\
u=\psi \text{ on } \Gamma_1,\\
g(u,n)=0, \ \ (\nabla_n u)^{\rm tan} + S_n(u) = w \text{ on } {\Gamma_2}.
\end{array}
\end{equation}
Only minor modifications need be made to our previous proofs, so we leave this
for the interested reader.

Of course, setting $\Gamma_2 = \emptyset$ proves the theorem when 
${\mathcal D}_{\mu,mix}^s$ is replaced by ${\mathcal D}_{\mu,D}^s$,
while setting $\Gamma_1 = \emptyset$ proves the theorem  in the case that
${\mathcal D}_{\mu,mix}^s$ is replaced by ${\mathcal D}_{\mu,N}^s$.

\section{The Stokes decomposition for manifolds with boundary}
\label{sec3}
In this section we recall well-known results about the Hodge decomposition
for manifolds with boundary (see \cite{Duff} and \cite{Mor} for proofs), and
define a new Stokes decomposition based on the solution to the Stokes problem,
whose summands are $\langle \cdot, \cdot \rangle_e$-orthogonal.

Let $(M,g)$ be a $C^\infty$ compact, oriented Riemannian $n$-dimensional manifold with
$C^\infty$ boundary ${\partial M}$, and let $i:{\partial M} \rightarrow M$ be
the inclusion map.  Then for a smooth vector field $X$ on $M$ and $n$, the
outward-pointing normal vector field on ${\partial M}$, 
$i^*(X\intprod\mu) = g(X,n) \mu_{\partial}$
where $\mu$ is the Riemannian volume form, and $\mu_{\partial}$ is the volume
form on ${\partial M}$ coming from the induced Riemannian metric.

By the trace theorem, $i^* \alpha$ is well-defined on ${\partial M}$ for 
$\alpha \in H^s(\Lambda^k(M))$ when $s \ge 1$; hence, for such $s$, 
$\alpha \in H^s(\Lambda^k(M))$ is {\it tangent} ($\parallel$) to ${\partial M}$
if and only if $n \intprod \alpha =0$, and
{\it normal} ($\perp$) to ${\partial M}$ if and only if $n^\sharp \wedge \alpha =0$.

When ${\partial M} = \emptyset$, $(d \alpha,\beta)_{L^2} = (\alpha, d \beta)_{L^2}$,
where $(\phi,\psi)_{L^2} = \int_M \phi\wedge *\psi$ 
(here, $*: \Lambda^k(M) \rightarrow \Lambda^{n-k}(M)$ 
denotes the Hodge star operator), and we have the standard Hodge decomposition
$$H^s(\Lambda^k) = d\bigl(H^{s+1}(\Lambda^{k-1}) \bigr) \oplus
\delta\bigl(H^{s+1}(\Lambda^{k+1}) \bigr) \oplus {\mathcal H}^{s,k},$$
where ${\mathcal H}^{s,k} = \{ \alpha \in H^s(\Lambda^k(M))| d \alpha=0 \text{ and }
\delta \alpha=0\}$ are the Harmonic fields.

When ${\partial M} \neq \emptyset$, we have that
$$(d \alpha, \beta)_{L^2} - (\alpha, \delta \beta)_{L^2} = \int_{\partial M}
(n^\sharp \wedge \alpha, \beta) \mu_{\partial}$$
and
$$(\delta \alpha, \beta)_{L^2} - (\alpha, d \beta)_{L^2} = -\int_{\partial M}
(n\intprod \alpha, \beta) \mu_{\partial}.$$
This shows that if $\delta \alpha =0$, then $\alpha \parallel {\partial M}$ iff
$(\alpha, d \beta)_{L^2} =0$ for all $\beta$, in which case the notion of 
$\alpha \parallel {\partial M}$ is well-defined even if $\alpha$ is only of class
$L^2$.  Similarly, if $d \alpha=0$, then $\alpha \perp {\partial M}$ iff
$(\alpha, d \beta)=0$ for all $\beta$.
We define 
\begin{align*}
H^s_t(\Lambda^k) &= \{\alpha \in H^s(\Lambda^k(M))\  |\ \alpha \parallel {\partial M}\},\\
H^s_n(\Lambda^k) &= \{\alpha \in H^s(\Lambda^k(M))\  |\ \alpha \perp {\partial M}\},\\
{\mathcal H}^{s,k}_t &= \{\alpha\in {\mathcal H}^s\  |\ \alpha \parallel {\partial M}\},\\
{\mathcal C}^{s,k}_t &= \{\alpha \in H^s(\Lambda^k(M))\  |\ \delta \alpha =0 \text{ and }
\alpha \parallel {\partial M}\}.
\end{align*}
Then for $s\ge 0$, we have the Hodge decompositions 
\begin{align*}
H^s(\Lambda^k) &= d\bigl(H_n^{s+1}(\Lambda^{k-1}) \bigr) \oplus
\delta\bigl(H_t^{s+1}(\Lambda^{k+1}) \bigr) \oplus {\mathcal H}^{s,k},\\
H^s(\Lambda^k) &= d\bigl(H^{s+1}(\Lambda^{k-1}) \bigr) \oplus {\mathcal C}^{s,k}_t,
\end{align*}
from which we can define the $L^2$ orthogonal projection onto ker$(\delta)$.

Consider the Hodge Laplacian $-\triangle= \delta d + d \delta$ with domain
$$\{ \alpha \in H^2(\Lambda^k(M)) \ | \ n \intprod \alpha =0 \text{ and }
n \intprod d \alpha =0 \},$$ 
and let ${\mathfrak P}_t$ denote the $L^2$ orthogonal
projection onto ${\mathcal H}^{s,k}_t$.  We call
\begin{equation}\nonumber
\begin{array}{c}
P_e: H^s(\Lambda^k) \rightarrow H^s_t(\Lambda^k)\\
P_e(\omega) = {\mathfrak P}_t \omega + \delta d (-\triangle)^{-1}(\omega - 
{\mathfrak P}_t\omega)
\end{array}
\end{equation}
the $L^2$ Hodge projection.

We shall now restrict our attention to $H^s(\Lambda^1(M))$ and identifying $1$-forms
with vector fields thru the metric $g$ on $M$.  Letting ${\mathcal X}^s_t
=\{u\in H^s(TM) \ | \ \text{div }u=0, v\parallel {\partial M} \}$, 
we may equivalently express the Hodge decomposition as 
$$H^s(TM) = \text{grad}H^{s+1}(M) \oplus {\mathcal X}^s_t,$$
so that for all $u\in H^s(TM)$, $u = v + \text{grad }p$, where $v\in{\mathcal X}^s_t$
and $p:M \rightarrow {\mathbb R}$ is obtained as the solution of Neumann problem
\begin{align*}
\triangle p &= \text{div } u  \quad \text{ in } M \\
g(\text{grad } p, n)& = g(u,n) \quad \text{ on } \partial M.
\end{align*}
Thus, a convenient and equivalent formula for the $L^2$ Hodge projection is
$$P_e:H^s(TM) \rightarrow {\mathcal X}^s_t, \ \ P_e(u) = u- \text{grad }p.$$

For each $\eta\in \mathcal{D}_\mu^s$, we define the projector
\begin{equation}\nonumber
\begin{array}{c}
P_\eta:T_\eta\mathcal{D}^s \rightarrow T_\eta \mathcal{D}_\mu^s, \\
P_\eta (X) = (P_e(X \circ \eta^{-1}))\circ \eta.
\end{array}
\end{equation}
Thus $\overline P: T{\mathcal D}^s \rightarrow T{\mathcal D}_\mu^s$, 
given on each fiber by $P_\eta$, is a bundle map covering the identity and is
$C^\infty$ by Appendix A of \cite{EM}.

Next, we define a new projector based on the elliptic Stokes problem.
Let  ${\mathfrak G}^s$ denote ${\mathcal D}_{D}^s$, 
${\mathcal D}_{N}^s$, or ${\mathcal D}_{mix}^s$, and similarly, let
${\mathfrak G}^s_\mu$ denote 
${\mathcal D}_{\mu,D}^s$, ${\mathcal D}_{\mu,N}^s$, or
${\mathcal D}_{\mu,mix}^s$. 

For $r\ge 1$, let ${\mathcal V}^r$ denote the $H^r$ vector fields on $M$
which satisfy the boundary conditions prescribed to elements of 
$T_e{\mathfrak G}^s$, and set ${\mathcal V}^r_\mu = \{ u \in {\mathcal V}^r \ |
\ \text{div }u=0\}$.   If $1\le r<2$, then elements of ${\mathcal V}^r$ and
${\mathcal V}^r_\mu$ only satisfy the essential boundary conditions
($u=0$ on ${\partial M}$ if ${\mathfrak G}_\mu^s={\mathcal D}_{\mu,D}^s$,
$g(u,n)=0$ on ${\partial M}$ if ${\mathfrak G}_\mu^s={\mathcal D}_{\mu,N}^s$,
or $u=0$ on $\Gamma_1$ and $g(u,n)=0$ on $\Gamma_2$
if ${\mathfrak G}_\mu^s={\mathcal D}_{\mu,mix}^s$) because vector fields in
${\mathcal V}^r$ for $r<2$ do not possess sufficient regularity for the
trace map to detect derivatives on the boundary.

We set ${\mathcal L}=-2\text{Def}^*\text{Def}$, and consider the
positive self-adjoint unbounded operator $(1-{\mathcal L})$ on $L^2(TM)$ 
with domain $D(1-{\mathcal L})={\mathcal V}^2$.  

\begin{prop}\label{P}
For $r\ge 1$ we have the following well defined decomposition
\begin{equation}\label{p1}
{\mathcal V}^r = {\mathcal V}^r_\mu \oplus (1-{\mathcal L})^{-1}{\rm grad}
H^{r-1}(M).
\end{equation}
Thus, if $F\in {\mathcal V}^r$, then there exists $(v,p) \in
{\mathcal V}^r_\mu \times H^{r-1}(M)/{\mathbb R}$ such that
$$F= v + (1-{\mathcal L})^{-1} \text{\rm grad }p $$
and the pair $(v,p)$ are solutions of the Stokes problem
\begin{equation}\label{p2}
\begin{array}{c}
(1-{\mathcal L}) v + \text{\rm grad }p = (1-{\mathcal L})F, \\
\text{\rm div }v =0,\\
v \text{ \rm satisfies boundary conditions}\\
\text{\rm prescribed to elements of } {\mathcal V}^r.\\
\end{array}
\end{equation}
The summands in (\ref{p1}) are $\langle \cdot, \cdot\rangle_e$-orthogonal.
Now, define the Stokes projector
\begin{equation}\label{p3}
\begin{array}{c}
{\mathcal P}_e: {\mathcal V}^r \rightarrow {\mathcal V}^r_\mu,\\
{\mathcal P}_e(F) = F - (1-{\mathcal L})^{-1}\text{\rm grad }p.
\end{array}
\end{equation}
Then, for $s>(n/2)+1$, $\overline {\mathcal P}: T{\mathfrak G}^s \rightarrow
T{\mathfrak G}^s_\mu$, given on each fiber by
\begin{equation}\nonumber
\begin{array}{c}
\overline{\mathcal P}_\eta: T_\eta{\mathfrak G}^s \rightarrow
T_\eta{\mathfrak G}^s_\mu,\\
\overline{{\mathcal P}}_\eta(X_\eta) = \left[ {\mathcal P}_e(X_\eta\circ
\eta^{-1})\right] \circ\eta,
\end{array}
\end{equation}
is a $C^\infty$ bundle map covering the identity.
\end{prop}
\begin{proof}
Acting on divergence-free vector-fields, ${\mathcal L} = \triangle_r$.
Thus, the proof that ${\mathcal D}_{\mu,mix}^s$ is a $C^\infty$ subgroup
of ${\mathcal D}_\mu^s$ shows that the
Stokes problem (\ref{p2}) has a unique solution $(v,p) \in
{\mathcal V}^r_\mu \times H^{r-1}(M)/{\mathbb R}$ for any $F\in {\mathcal V}^r$,
$r \ge 1$.

It is easy to verify that the summands in (\ref{p1}) are $\langle\cdot,\cdot
\rangle_e$-orthogonal, so it only remains to show that $\overline {\mathcal P}$
is smooth.   For $F_\eta \in T_\eta {\mathfrak G}^s$, let $F
=F_\eta\circ\eta^{-1}$, and let $(v,p)$ solve (\ref{p2}).  
By (\ref{p3}), it suffices to prove that 
$$\left[(1-{\mathcal L})^{-1}\text{grad }p\right] \circ \eta =
\left[
(1-{\mathcal L})^{-1} \text{grad}\triangle^{-1} \text{div}(1-{\mathcal L})
(v-F) \right] \circ \eta
$$
is smooth.  Letting $V_\eta = v \circ\eta \in T_\eta{\mathfrak G}^s_\mu$,
we have the equivalent expression for
$\left[(1-{\mathcal L})^{-1}\text{grad }p\right] \circ \eta$ given by
$$
\overline{(1-{\mathcal L})^{-1}}_\eta \circ \overline{\text{grad}}_\eta \circ
\overline{\triangle^{-1}}_\eta \circ \overline{\text{div}}_\eta \circ
\overline{(1-{\mathcal L})}_\eta (V_\eta - F_\eta)
$$
which is a $C^\infty$ bundle map by Proposition \ref{smooth_bmap} together
with Lemmas \ref{b1} and \ref{b1a}.
\end{proof}

\section{A new right invariant metric on ${\mathcal D}_{\mu,D}^s$,
${\mathcal D}_{\mu,N}^s$, ${\mathcal D}_{\mu,mix}^s$ and its geodesics}
\label{sec4}

Recall that a {\it weak Riemannian metric} on a Hilbert manifold ${\mathcal M}$
is given by a map $\gamma$ which assigns to each $m \in {\mathcal M}$, a 
continuous positive-definite symmetric bilinear form $\gamma(m) \in T^*_m{\mathcal M}
\otimes T^*_m{\mathcal M}$, which is $C^\infty$ with respect to $m\in {\mathcal M}$.
The metric $\gamma$ is termed {\it weak}, because it defines a topology which is
weaker than the original topology on ${\mathcal M}$ (and hence on 
$T_m{\mathcal M}$).

In general, the geodesic flow of a weak metric does not exist. A simple example
is given by the lack of a well-defined exponential map for the usual $L^2$ metric 
on ${\mathcal D}^s$ when ${\partial M}$ is not empty.
Nevertheless,
the seminal paper of Ebin-Marsden \cite{EM} proves that 
it is indeed possible to define a weak right invariant $L^2$ metric on 
${\mathcal D}_\mu^s$ for manifolds with boundary, and that this weak metric 
induces a (weak) Levi-Civita covariant derivative and geodesic
flow.  As we have described, the geodesic flow of the invariant $L^2$ metric on
${\mathcal D}_\mu^s$ generates solutions to the Euler equations of ideal
hydrodynamics; we shall introduce a new weak invariant metric on ${\mathcal D}_\mu^s$
which, remarkably, also generates geodesic flow that solves the equations of
ideal non-Newtonian second-grade fluids as well as the averaged Euler or Euler-$\alpha$
equations.

Let ${\mathfrak G}^s_\mu$ denote either ${\mathcal D}_{\mu,D}^s$,
${\mathcal D}_{\mu,N}^s$, or ${\mathcal D}_{\mu,mix}^s$, and let  
$\bar g$ denote the induced inner-product on the fibers of 
$[T^*M$$\otimes$ ${T^*M]^*}^{\otimes 2}$. 

\begin{prop}
Define the bilinear form $\langle\cdot, \cdot\rangle_e$ on 
$T_e {\mathfrak G}^s_\mu$ as follows:
for $X,Y \in T_e {\mathfrak G}^s_\mu$ and $\alpha >0$, set
\begin{align}
\langle X, Y \rangle_e = & \int_M \bigl(g_x(X(x),Y(x))  +\frac{\alpha^2}{2} 
\bar g_x (\pounds_X g(x), \pounds_Y g(x)) \bigr)\mu(x),
\label{metric}
\end{align}
and define a bilinear form on each fiber of $T{\mathfrak G}^s_\mu$ by
right translation so that for $X_\eta, Y_\eta
\in T_\eta {\mathfrak G}^s_\mu$,
$$ \langle X_\eta, Y_\eta \rangle_\eta = \langle X_\eta \circ \eta^{-1}, Y_\eta \circ
\eta^{-1} \rangle_e.$$
Then $\langle \cdot, \cdot \rangle$, given on each fiber by $\langle \cdot, 
\cdot \rangle_\eta$, is a right invariant weak Riemannian metric on ${\mathfrak G}^s_\mu$.
\end{prop}
\begin{proof}
That $\langle \cdot, \cdot \rangle$ is $C^\infty$ on 
${\mathfrak G}^s_\mu$ follows from Lemma \ref{o-lemma}.  
That $\langle \cdot, \cdot
\rangle_\eta$ is a positive-definite symmetric bilinear form is proven as follows:
$$ 2\text{Def }u = \pounds_u g = \nabla u +  \nabla u^t,$$
so for any of the boundary conditions prescribed
on elements of $T_e{\mathfrak G}^s_\mu$, we have that
$$ 0 \le 2\text{Def}^* \text{Def }u = -(\triangle + 2 \text{Ric})u,$$
so that integrating by parts (and noting that the boundary terms vanish), we may
express $\langle \cdot, \cdot \rangle_e$ in the equivalent form
$$
\langle X, Y \rangle_e =  \int_M g_x\bigl((1 - \triangle_r) X(x),Y(x)\bigr) \mu(x).
$$
Since $(1-\triangle_r)$ is a self-adjoint positive operator (on $L^2$ vector
fields that are divergence-free), this shows that
$\langle \cdot, \cdot \rangle$ is a well defined $C^\infty$ weak invariant 
Riemannian metric on ${\mathfrak G}_\mu^s$.
\end{proof}
The metric $\langle \cdot, \cdot \rangle$ is invariant under the action of
${\mathfrak G}^s_\mu$, so the subgroups of the  volume preserving diffeomorphism group
that we have constructed play the
role of both configuration space as well as symmetry group (this is the massive
particle relabeling symmetry of hydrodynamics).  In order to formally establish
the equations of geodesic motion of the invariant metric $\langle \cdot, \cdot \rangle$
on ${\mathfrak G}^s_\mu$ we shall make use of the Euler-Poincar\'{e} reduction theorem.
The reader unfamiliar with this symmetry reduction procedure is referred to Appendix
\ref{app1} for a brief discussion.

\begin{prop}\label{EP}
Let the pair $({\mathfrak G}^s_\mu, \langle \cdot , \cdot \rangle)$ denote either
${\mathcal D}_{\mu,D}^s$, ${\mathcal D}_{\mu,N}^s$, or
${\mathcal D}_{\mu,mix}^s$ together with the right invariant Riemannian metric
defined in (\ref{metric}).  Then, a curve $\dot\eta(t) \in T{\mathfrak G}^s_\mu$ is
a geodesic of $\langle \cdot, \cdot \rangle$ if and only if  its projection onto
the fiber over the identity given by $u(t) = \dot\eta(t) \circ \eta(t)^{-1} \in
T_e{\mathfrak G}^s_\mu$ is a solution of
\begin{equation}\label{ep}
\begin{array}{c}
(1-\alpha^2\triangle_r) \partial_t u + \nabla_u(1-\alpha^2 \triangle_r) u 
- \alpha^2 \nabla u^t \cdot \triangle_r u = -\text{\rm grad } p, \\
\text{\rm div }u =0, \ \ \ u(0)=u_0,
\end{array}
\end{equation}
together with the boundary conditions
\begin{align*}
u=0 \text{ \rm on }{\partial M} & \ \ \text{ \rm if } \ \
           {\mathfrak G}^s_\mu = {\mathcal D}_{\mu,D}^s, \\
\\
\left. \begin{array}{cc}
g(u,n)=0,\\
(\nabla_n u)^{\rm tan} + S_n(u)=0 
\end{array}\right\}
\text{ \rm on }{\partial M} & \ \ \text{ \rm if } \ \
           {\mathfrak G}^s_\mu = {\mathcal D}_{\mu,N}^s, \\
\\
\left.  \begin{array}{cl}
                u=0 & \text{ \rm on } \Gamma_1 \\
                \left. \begin{array}{cc}
                g(u,n)=0,\\
                (\nabla_n u)^{\rm tan} + S_n(u)=0 
                \end{array} \right\}
                & \text{ \rm on } \Gamma_2
                \end{array}
       \right\} & \ \ \text{ \rm if }  \ \ {\mathfrak G}^s_\mu = 
                                        {\mathcal D}_{\mu,mix}^s,
\end{align*}
where $\text{\rm grad }p$ is completely determined by the Stokes projector
${\mathcal P}_e$.
\end{prop}
\begin{proof}
From part (d) of Proposition \ref{thm_ep}, the reduced Lagrangian is given by
$\langle \cdot , \cdot \rangle_e$, so that $\dot\eta(t)$ is a geodesic of
$\langle \cdot, \cdot \rangle$ on ${\mathfrak G}^s_\mu$ if $u(t) =
\dot\eta(t) \circ \eta(t)^{-1}$ is a fixed point of the reduced action function
(on an arbitrary interval $(a,b)$)
$s:T_e{\mathfrak G}^s_\mu \rightarrow {\mathbb R}$ given by
$$ s(u) = {\frac{1}{2}}\int_a^b \langle u(t),u(t) \rangle_e dt.$$

Let $\epsilon \mapsto \eta^\epsilon$ be a smooth curve in ${\mathfrak G}^s_\mu$
such that $\eta^0=\eta$ and $(d/d \epsilon)\eta^\epsilon|_{\epsilon=0} 
= \delta \eta \in
T_\eta {\mathfrak G}^s_\mu$; the map $t \mapsto \delta \eta(t)$ is the 
variation of the curve $\eta(t)$ on the interval $(a,b)$ and 
$\delta \eta(a)= \delta\eta(b)=0$.  The curve $\epsilon \mapsto \eta^\epsilon$
induces a curve $\epsilon \mapsto u^\epsilon$ in the single fiber
$T_e{\mathfrak G}^s_\mu$ such that $u^0=u$ and 
$(d/d \epsilon)u^\epsilon|_{\epsilon=0} = \delta u$.  The Euler-Poincar\'{e}
reduction theorem gives the relation
$$\delta u = \partial_t(\delta \eta \circ \eta^{-1}) + [ \delta \eta \circ
\eta^{-1}, u]_e.$$

Computing the first variation of the action $s$, we have that
\begin{align*}
ds(u)&\cdot \delta u  \\
=&\int_a^b \int_M  \left(
g(u,\delta u)+2\alpha^2\bar g(\text{Def }u,\text{Def }\delta u) \right)\mu d t\\
=&\int_a^b \left[\int_M  g((1-\alpha^2\triangle_r)u, \delta u) \mu
   + \alpha^2\int_{\partial M}  g( (\nabla_nu)^{\rm tan} + S_n(u), \delta u) 
       \mu_\partial \right] dt.
\end{align*}
Since $u$ and $\delta u$ satisfy the boundary conditions prescribed to
elements of $T_e{\mathfrak G}^s_\mu$, 
the boundary term in the above equation vanishes, leaving only
$$ ds(u) \cdot \delta u =
\int_a^b \int_M  g\bigl((1-\alpha^2\triangle_r)u, 
\partial_t(\delta \eta \circ \eta^{-1}) + [ \delta \eta \circ
\eta^{-1}, u]_e\bigr) \mu dt.$$
Using the formula $[x,y]_e= \nabla_yx - \nabla_xy$ and integrating by parts,
we obtain
\begin{align*}
ds(u)& \cdot \delta u  \\
=&\int_a^b \int_M  g\bigl( (1-\alpha^2\triangle_r)\partial_t u +
\nabla_u(1-\alpha^2\triangle_r)u - \alpha^2\nabla u^t\cdot \triangle_r u, 
\delta \eta \circ \eta^{-1} \bigr) \mu dt \\
=& \int_a^b \bigl\langle \partial_t u + (1-\alpha^2{\mathcal L})^{-1} \bigl[
\nabla_u(1-\alpha^2\triangle_r)u - \alpha^2\nabla u^t \cdot \triangle_r u\bigr],
\delta \eta \circ \eta^{-1} \bigr\rangle_e dt,
\end{align*}
where again ${\mathcal L} = -2\text{Def}^*\text{Def}$.
Since right translation is an isomorphism, $\delta \eta \circ \eta^{-1}
\in T_e {\mathfrak G}^s_\mu$ is arbitrary, so $u$ is a fixed point of $s$ iff
$$\partial_t u + {\mathcal P}_e \left(
(1-\alpha^2{\mathcal L})^{-1} \bigl[
\nabla_u(1-\alpha^2\triangle)u - \alpha^2\nabla u^t \cdot \triangle u \bigr]
\right) =0,$$
and this is precisely (\ref{ep}), as $(1-\alpha^2{\mathcal L}) \partial_t u
=(1-\alpha^2 \triangle_r) \partial_t u$ since $\text{div}\partial_t u =0$.
\end{proof}

In the next section, we prove Theorem \ref{Thm2} by establishing existence 
and uniqueness of geodesics of the invariant metric. 
The following simple lemma will play a fundamental role.

\begin{lemma}\label{commute}
For $s>(n/2)+1$, let $u,v \in T_e {\mathcal D}_{\mu,D}^s$, and
consider the unbounded self-adjoint operator $(1-{\mathcal L})$ on $L^2$
with domain $D(1-{\mathcal L}) = H^2(TM) \cap H^1_0(TM)$.
Then
\begin{align*}
(1-{\mathcal L}) \nabla_u v = & \nabla_u (1-\triangle_r)v - 
{\rm div}[\nabla v \cdot \nabla u^t + \nabla v \cdot \nabla u] 
-\text{\rm grad Tr}[\nabla u \cdot \nabla v] \\
& + (\nabla_u {\rm Ric}) \cdot v -\text{\rm grad Ric}(u,v)
- {\rm Tr}[ \nabla (R(u,\cdot)v)+R(u,\cdot)\nabla v].
\end{align*}
\end{lemma}
\begin{proof}
First notice that for $s>(n/2)+1$, $\nabla_uv$ is an $H^{s-1}$ vector field
on $M$ whose trace vanishes on ${\partial M}$; thus, it makes sense for
the operator $(1-{\mathcal L})$ to act on $\nabla_uv$.

Recall that ${\mathcal L} = -(\triangle + 2\text{Ric} + \text{grad div})$,
so we begin by computing the commutator of $[-\triangle, \nabla_u]$.  
Let $\{e_i\}$ be a local orthonormal frame, and write the Hodge Laplacian 
$\triangle =-( d \delta + \delta d)$ acting on $1$-forms 
(identified with vector fields) as $\triangle = \nabla_{e_i}\nabla_{e_i} +
\text{Ric}$, so that
$$ 
\triangle \nabla_u v = \nabla_{e_i} \nabla_{e_i} (\nabla_u v) - \text{Ric}(\nabla_uv).
$$
Using the definition of the Riemannian curvature operator, we compute that
\begin{align*}
\nabla_{e_i} \nabla_{e_i} \nabla_u v  & =  \nabla_{e_i} \left( R(u,e_i)v\right) + 
\nabla_{e_i}( \nabla_{[e_i,u]}v ) + \nabla_{e_i} \nabla_u \nabla_{e_i} v\\
&= \nabla_u \nabla_{e_i} \nabla_{e_i} v + \nabla_{e_i} ( \nabla_{[e_i,u]}v)
+\nabla_{[e_i,u]} \nabla_{e_i}v \\
&\qquad + \nabla_{e_i}\left( R(u,e_i)v \right) + R(u,e_i) \nabla_{e_i} v +
\nabla_u \text{Ric}(v) - \nabla_u\left( \text{Ric}(v)\right).
\end{align*}
Expressing $u$ as $u^je_j$, we see that $[e_i,u] = e_i[u^j]e_j$; hence, one may
easily verify that
\begin{align*}
\nabla_{[e_i,u]} \nabla_{e_i} v & = \text{div}[ \nabla v \cdot \nabla u],\\
\nabla_{e_i} \left(\nabla_{[e_i,u]} v \right)& = \text{div}[ \nabla v \cdot \nabla u^t],
\end{align*}
so that
\begin{align*}
-\triangle \nabla_u v =& -\nabla_u \triangle v -
\text{div}\left[ \nabla v \cdot \nabla u^t + \nabla v \cdot \nabla u\right]
-(\nabla_u\text{Ric}) \cdot v \\
& \qquad - \text{Tr}\left[ \nabla(R(u,\cdot) v) + R(u\cdot)\nabla v \right].
\end{align*}

Using the fact that div$\nabla_uv = \text{Tr}( \nabla u \cdot \nabla v) +
\text{Ric}(u,v)$, and combining terms involving the Ricci curvature gives the
result.
\end{proof}

We remark that if we embed $M$ into its double $\tilde M$, smoothly extending 
$g$, and let $(1-\hat{\mathcal L})$ denote the operator 
$(1-2\text{Def}^*\text{Def})$ on $\tilde M$, then it makes sense for
$R \circ (1-\hat{\mathcal L}) \circ E$ to formally act on an arbitrary vector 
fields
on $\overline M$. Here,  $R$ denotes restriction and $E$ denotes extension;
see the proof of Theorem \ref{Thm2} for a more detailed construction of
such an operator. It follows that the above lemma also holds for the groups
${\mathcal D}_{\mu,N}^s$ and ${\mathcal D}_{\mu,mix}^s$  when the operator
$(1-{\mathcal L})$ acting vector fields which vanish on ${\partial M}$ is
replaced by $R \circ (1-\hat{\mathcal L})\circ E$.  

\section{Proof of Theorem \ref{Thm2}}
\label{sec5}
Let us denote the covariant material time derivative by $(\nabla / dt)$. For the
remainder of this section we shall, for convenience, set $\alpha =1$.  The
unbounded, self-adjoint operator $(1-{\mathcal L})=(1-2\text{Def}^*\text{Def})$ 
on $L^2(TM)$ has domain $H^2(TM)\cap H^1_0(TM)$.

\begin{prop}\label{geodesic}
For $s>(n/2)+1$, let $\eta(t)$ be a curve in ${\mathcal D}_{\mu,D}^s$,
and set $u(t)=\dot\eta \circ
\eta(t)^{-1}$.  Then $u$ is a solution of the initial-boundary value problem
(\ref{avg_euler}) with Dirichlet boundary conditions $u=0$ on ${\partial M}$ 
if and only if
\begin{equation}\label{simple}
\overline{\mathcal P}_\eta \circ \left[
\frac{\nabla \dot \eta}{dt} + \left[-\nu(1-{\mathcal L})^{-1}\triangle_r u +
{\mathcal U}(u) + {\mathcal R}(u)\right]
\circ \eta \right] =0,
\end{equation}
where
\begin{align*}
{\mathcal U}(u) =& (1-{\mathcal L})^{-1}\bigl\{
\text{\rm div}\left[ \nabla u \cdot 
\nabla u^t + \nabla u \cdot \nabla u - \nabla u^t \cdot \nabla u\right] 
+\text{\rm grad Tr}(\nabla u \cdot \nabla u)\bigr\} \\
{\mathcal R}(u) =& (1-{\mathcal L})^{-1} \bigr\{  \text{\rm Tr}
\left[ \nabla \left( R(u,\cdot)u \right) +R(u,\cdot) \nabla u + 
R(\nabla u, \cdot)u \right] \\
&\qquad\qquad\qquad  + \text{\rm grad Ric}(u,u)
-(\nabla_u\text{\rm Ric}) \cdot u + \nabla u^t \cdot \text{\rm Ric}(u)
\bigr\}  ,
\end{align*}
and $\overline{\mathcal P}_\eta: T_\eta{\mathcal D}_D^s \rightarrow 
T_\eta{\mathcal D}_{\mu,D}^s$ is the Stokes projector.
\end{prop}
\begin{proof}
We first set $\nu=0$.
Covariantly differentiating ${\dot{\eta}} = u \circ \eta$ yields
$$\frac{\nabla {\dot{\eta}}}{dt} \circ \eta^{-1} = \partial_t u + \nabla_u u.$$
Using Lemma \ref{commute}, we obtain that
\begin{align*}
(1-{\mathcal L}) \left(\frac{\nabla \dot{\eta}}{dt} \circ \eta^{-1} \right) = &
(1-\triangle_r) \partial_t u +  (1-{\mathcal L})\nabla_u u \\
=& (1-\triangle_r) \partial_t u + \nabla_u(1-\triangle_r)u -\text{div}[\nabla u \cdot
\nabla u^t + \nabla u \cdot \nabla u] \\
&-\text{grad Tr}(\nabla u \cdot\nabla u) - \text{grad Ric}(u,u) \\
& - \text{Tr} \left[ \nabla( R(u,\cdot)u) + R(u,\cdot)\nabla u\right] + \left(
\nabla _u \text{Ric}\right) \cdot u.
\end{align*}
Now $\nabla u^t \cdot \triangle u =\text{div}[ \nabla u^t \cdot \nabla u] +
\text{grad } \phi - \text{Tr}R(\nabla u, \cdot) u - \nabla u^t \cdot \text{Ric}(u)$,
for some $\phi: M \rightarrow {\mathbb R}$; hence,
$$(1-\triangle_r) \partial_t u + \nabla_u(1-\triangle_r)u - 
\nabla u^t \cdot\triangle u = -\text{grad }p$$
if and only if
$$\frac{\nabla {\dot{\eta}}}{dt} \circ \eta^{-1} + {\mathcal U}(u) +
{\mathcal R}(u) = -(1-{\mathcal L})^{-1}\text{grad }\tilde p,$$
for some $\tilde p : M \rightarrow {\mathbb R}$,
and by Proposition \ref{P}, this is precisely equation (\ref{simple}) 
with $\nu=0$.   Adding the term $ {\mathcal P}_\eta \circ 
-\nu(1-{\mathcal L})^{-1}\triangle_r u \circ \eta$ 
to (\ref{simple}) produces the equation (\ref{avg_euler}).
\end{proof}

We can now proceed with the proof of the theorem.  We first consider
the inviscid case first with the viscosity $\nu=0$.

By Proposition \ref{geodesic}, the geodesic flow of the invariant metric
$\langle \cdot, \cdot \rangle$ is the solution of
$$\frac{\nabla \dot\eta}{dt} =\overline{\mathcal S}_\eta(\dot\eta):=
(1-{\mathcal P}_\eta)\frac{\nabla\dot\eta}{dt} - {\mathcal P}_\eta\circ
( \overline {\mathcal U}_\eta + \overline{\mathcal R}_\eta) \dot\eta,$$
where  $\overline {\mathcal S}$ is the bundle map covering the identity
given on each fiber by $\overline {\mathcal S}_\eta$, and
$$ \overline{\mathcal U}_\eta(X_\eta) = [ {\mathcal U}(X_\eta \circ
\eta^{-1})] \circ \eta, \ \
\overline{\mathcal R}_\eta(X_\eta) = [ {\mathcal R}(X_\eta \circ
\eta^{-1})] \circ \eta \ \ \forall \ X_\eta \in T_\eta{\mathfrak G}^s_\mu.$$

Now the second tangent bundle $T^2 {\mathcal D}_{\mu,D}^s$ is identified
with $H^s$ maps ${\mathcal Y}: M \rightarrow T^2M$ which cover some
$X_\eta\in T_\eta {\mathcal D}_{\mu,D}^s$.  The second-order vector field
$\ddot \eta: M \rightarrow T^2M$ is just such a
map, covering $\dot\eta\in T_\eta{\mathcal D}_{\mu,D}^s$.

Using a local representation, 
we may express the material time derivative above as the system
\begin{align*}
\dot\eta &= V_\eta, \\
\ddot \eta &= \frac{dV_\eta}{dt} = {\mathcal B}(\eta,\dot\eta) =
-\Gamma_\eta(\dot\eta,\dot\eta) + \overline{\mathcal S}_\eta(\dot\eta), \\
\eta&(0)=e,\\
V_\eta&(0)=u_0,\\
\end{align*}
since $\nabla {\dot{\eta}}/dt = \ddot \eta + 
\Gamma_{\eta}({\dot{\eta}},{\dot{\eta}})$, where $\Gamma_\eta(\dot\eta,
\dot\eta)$ is 
the Christoffel map, given in a local coordinate chart on $M$ by
$\Gamma_{\eta(x)}({\dot{\eta}},{\dot{\eta}}) = \Gamma^i_{jk}(x)
(\dot\eta \circ \eta^{-1})^j(\dot\eta \circ \eta^{-1})^k$.
${\mathcal B}(\eta,\dot\eta)$ is the principal part of the geodesic
spray of $\langle \cdot, \cdot\rangle$ on ${\mathcal D}_{\mu,D}^s$;
hence, with ${\mathcal U}$ denoting a local open neighborhood of
$\eta\in {\mathcal D}_{\mu,D}^s$, to establish the first assertion we shall 
prove that 
${\mathcal B}$ maps ${\mathcal U} \times H^s_\eta(TM)$ into
$H^s_\eta(TM)$, and that ${\mathcal B}$ is $C^\infty$.
The result then follows by application of the fundamental theorem of ordinary 
differential equations on Hilbert manifolds (see \cite{La}, Theorem 2.6),
and the existing time-reversal symmetry  $t \mapsto -t$.

As the Christoffel map is a $C^\infty$ map of ${\mathcal U} \times
H^s_\eta(TM)$ into $H^s_\eta(TM)$ (since $g$ is $C^\infty$ and $H^s$ 
is a multiplicative algebra), we must show that 
$\overline{\mathcal S}_\eta$ is $C^\infty$.
Since  $\overline{\mathcal P}_\eta: T_\eta{\mathcal D}_D^s \rightarrow 
T_\eta{\mathcal D}^s_{\mu,D}$ is $C^\infty$ by Proposition \ref{P},
to show that ${\mathcal P}_\eta \circ \overline {\mathcal U}_\eta :
T_\eta {\mathcal D}^s_{\mu,D} \rightarrow T_\eta {\mathcal D}^s_{\mu,D}$ is 
$C^\infty$ it suffices to prove that 
$$\overline{(1-{\mathcal L})^{-1}}_\eta \circ 
\overline{{\rm div}}_\eta \circ [\nabla (\dot{\eta} \circ \eta^{-1})\circ \eta \
 \cdot \ \nabla (\dot{\eta} \circ \eta^{-1})\circ \eta ]
: T_\eta {\mathcal D}^s_{\mu,D} \rightarrow T_\eta {\mathcal D}^s_D$$
and
$$\overline{(1-{\mathcal L})^{-1}}_\eta \circ 
\overline{{\rm grad}}_\eta \circ \text{Tr}[\nabla (\dot{\eta} \circ 
\eta^{-1})\circ \eta \
 \cdot \ \nabla (\dot{\eta} \circ \eta^{-1})\circ \eta ]
: T_\eta {\mathcal D}^s_{\mu,D} \rightarrow T_\eta {\mathcal D}^s_D$$
are $C^\infty$ bundle maps.
But this follows from Lemmas \ref{b1} and \ref{b1a} together with Proposition
\ref{smooth_bmap}.
Since $R$ and Ric are
$C^\infty$ on $M$, a similar argument shows that ${\mathcal P}_\eta \circ
\overline{\mathcal R}_\eta:
T_\eta {\mathcal D}^s_{\mu,D} \rightarrow T_\eta {\mathcal D}^s_{\mu,D}$ 
is $C^\infty$ as well.

We next prove that $(1-\overline {\mathcal P}_\eta) \circ (\nabla \dot\eta/dt)$
is $C^\infty$.  Since $\partial_t u \in T_e {\mathcal D}_{\mu,D}^s$,
$$ \overline{\mathcal P}_\eta\circ \frac{\nabla \dot\eta}{dt}
= \left[ \partial_t u + {\mathcal P}_e(\nabla_uu)\right] \circ \eta,$$
so that 
$$(1-\overline{\mathcal P}_\eta) \circ (\nabla \dot\eta/dt)
= -(1-{\mathcal L})^{-1} \text{grad }p \circ \eta,$$
where $p$ depends on $v$ and the pair $(v,p)$ is a solution of the Stokes
problem
\begin{equation}\nonumber
\begin{array}{c}
(1-\triangle_r) v + \text{grad } p = (1-{\mathcal L}) \nabla_uu \\
\text{div } v=0\\
v=0 \text{ on } {\partial M}.
\end{array}
\end{equation}

Since $s>(n/2)+1$, $(1-{\mathcal L}) \nabla_uu$ is in $H^{s-3}(TM)$; the
argument in Step 3 of the proof of Theorem \ref{Thm2} then gives
a unique solution $(v,p)\in {\mathcal V}^{s-1}_\mu \times H^{s-2}(M)/
{\mathbb R}$.  If $-1< s-3 < 0$, then the pair $(v,p)$ is 
interpreted as a weak solution.  

A priori, $(1-{\mathcal L})^{-1}\text{grad } p$ is only in $H^{s-1}$, but
we shall show that, in fact, 
$(1-{\mathcal L})^{-1}\text{grad } p$ is actually of class $H^s$.
We have that
$$(1-{\mathcal L})^{-1}\text{grad }p = 
(1-{\mathcal L})^{-1}\text{grad}\triangle^{-1}\text{div}
(1-{\mathcal L})(v -\nabla_uu).$$
We embed $M$ into its double $\tilde M$, extending $g$ to $\tilde M$, and
choose a $C^\infty$ extension of $u$ to $\tilde M$.  
For any vector bundle ${\mathcal E}$ over $M$, let
$$E: H^s({\mathcal E} \downarrow M) \rightarrow H^s({\mathcal E} \downarrow
\tilde M), \ \ \ E(\xi)|_M = \xi$$
denote the linear extension operator, and let $R$ denote the corresponding
restriction operator.  Let $\hat {\mathcal L}$ denote $R \circ {\mathcal L}
\circ E$; then it makes sense to form the commutator of the operators
div with $\hat {\mathcal L}$, and the operator
$$[\text{div},\hat {\mathcal L}] : H^r(TM) \rightarrow H^{r-2}(TM)$$
is continuous.  Notice that as ${\mathcal L}$ is a local operator, if
$w=0$ on $M$, then $\hat {\mathcal L} w =0$ by the property of the extension
operator given above.  Since $\text{div }v=0$, 
$$-(1-{\mathcal L})^{-1}\text{grad}\triangle^{-1} \text{div}{\mathcal L}v
= -(1-{\mathcal L})^{-1}\text{grad}\triangle^{-1} [\text{div},\hat{\mathcal L}]
v,$$
which is in $H^s(TM)\cap H^1_0(TM)$, since
$$-(1-{\mathcal L})^{-1}\text{grad}\triangle^{-1} [\text{div},\hat{\mathcal L}]:
H^{s-1}(TM) \cap H^1_0(TM) \rightarrow H^s(TM) \cap H^1_0(TM)$$ 
is a compact operator.

The identical argument shows that
$-(1-{\mathcal L})^{-1}\text{grad}\triangle^{-1} [\text{div},\hat{\mathcal L}]
\nabla_uu$ is in $H^s(TM) \cap H^1_0(TM)$, since $\nabla_uu$ is in $H^{s-1}(TM)
\cap H^1_0(TM)$.  Since $\text{div}\nabla_uu =  
\text{Tr}(\nabla u \cdot \nabla u) + \text{Ric}(u,u)$ is an $H^{s-1}$ vector
field on $M$, and since
$$-(1-{\mathcal L})^{-1}\text{grad}\triangle^{-1}\hat{\mathcal L}: H^{s-1}(TM)
\rightarrow H^s(TM) \cap H^1_0(TM)$$
compactly, we see that
$$-(1-{\mathcal L})^{-1}\text{grad}\triangle^{-1} \text{div}{\mathcal L}
\nabla_uu$$
is in fact of class $H^s$.  Regularity up to the boundary immediately
follows from the fact that $\nabla_uu=0$ on $\partial M$.
Thus $(1-{\mathcal L})^{-1}\text{grad }p$ is in
$H^s$, and from Section \ref{app2}, it follows that
$[(1-{\mathcal L})^{-1}\text{grad }p] \circ \eta$ is in $H^s_\eta(TM)$.

The fact that $u$
is the unique solution of (\ref{avg_euler}) with $\nu=0$
is the statement of Proposition \ref{geodesic}.  That
$u$ is in $C^0( I, {\mathcal V}^s_\mu) \cap C^1(I,{\mathcal V}^{s-1}_\mu)$
and depends continuously on the initial data $u_0$ follows from
the fact that the inversion map
$(\eta \mapsto \eta^{-1}): {\mathcal D}^s \rightarrow {\mathcal D}^s$ is only
$C^0$, but is $C^1$ when considered as a map from ${\mathcal D}^s$ into
${\mathcal D}^{s-1}$.

This proves the theorem for the case $\nu=0$.

Next, we consider the viscous case $\nu>0$.
We need only show that the viscous term, thought of a bundle map,
$\dot\eta\mapsto \overline{\mathcal P}_\eta [(1-{\mathcal L})^{-1} \triangle_r
(\dot\eta\circ \eta^{-1})] \circ \eta : 
T_\eta{\mathcal D}_{\mu,D}^s \rightarrow T_\eta{\mathcal D}_{\mu,D}^s$
is a $C^\infty$ bundle map.  But this map is the same as 
$\overline{\mathcal P}_\eta \circ \overline{(1-{\mathcal L})^{-1}}_\eta
\circ (\overline{\triangle_r})_\eta$, which is a $C^\infty$ bundle map by
Proposition \ref{P} and Proposition \ref{smooth_bmap}.

The viscosity destroys the time-reversal symmetry, so the solution is
now defined on $\bar I$.  This concludes the proof of Theorem \ref{Thm2}.

\section{Proof of Theorem \ref{Thm3}}
\label{sec6}

The existence of the unique Levi-Civita covariant derivative of the
right invariant metric $\langle \cdot, \cdot \rangle$ on 
${\mathcal D}_{\mu,D}^s$ is an immediate consequence of the smoothness of
the geodesic flow of $\langle \cdot, \cdot \rangle$ provided by Theorem
\ref{Thm2}.  The formulas for $\tilde \nabla$ then follow from the fundamental
theorem of Riemannian geometry.

As to the properties of the curvature operator,
right invariance of $\tilde R$ follows from the right invariance of $\tilde 
\nabla$.  Next we prove that $\tilde R$ is bounded in $H^s$ for $s>
\frac{n}{2}+2$.  

Extend $X_\eta, Y_\eta, Z_\eta \in T_\eta{\mathcal D}^s_{\mu,D}$ to smooth 
right invariant 
vector fields $x^r, y^r,z^r$ on ${\mathcal D}^s_{\mu,D}$ and let $x=x^r(e),
y=y^r(e)$, and $z=z^r(e)$.  Let 
\begin{align*}
M_xy =& (1-{\mathcal P}_e)\circ \nabla_xy + 
(1/2) {\mathcal P}_e \circ ({\mathfrak U}(x,y) + {\mathfrak R}(x,y)).
\end{align*}
As the proof of Theorem \ref{Thm2} shows, $M$ has the following property:  
\begin{itemize}
\item[]
If $x$ and $y$ are $H^s$ divergence-free vector fields on $M$, and $s$ is 
sufficiently large so that $H^{s-1}(TM)$ forms a multiplicative algebra, then
there exists a positive constant $c$, such that
$|M_xy|_s \le c|x|_s|y|_s$.
\end{itemize}

Now, since $\tilde \nabla$ is right invariant, we have that
\begin{align*}
\tilde R_\eta&(X_\eta,Y_\eta)z^r_\eta = \left( \tilde \nabla_{y^r} \tilde \nabla_{x^r} z^r
\right)_\eta - \left( \tilde \nabla_{x^r} \tilde \nabla_{y^r} z^r \right)_\eta
+ \left( \tilde \nabla_{[x^r,y^r]} z^r \right)_\eta \\
&=\left[ (\nabla_y + M_y)(\nabla_x + M_x)z\right] \circ \eta
- \left[ (\nabla_x + M_x)(\nabla_y + M_y)z\right] \circ \eta \\
&\qquad \qquad + \left[ (\nabla_{[x,y]} + M_{[x,y]}) z\right] \circ \eta \\
&=\left[ (\nabla_y\nabla_x - \nabla_x\nabla_y + \nabla_{[x,y]})z\right]\circ 
\eta
+\left[ (M_yM_x - M_xM_y + M_{[x,y]})z\right]\circ \eta\\
&\qquad \qquad
+\left[ \{ \nabla_x,M_y\} z + \{M_x,\nabla_y\}z\right] \circ \eta ,
\end{align*}
where $\{ \cdot, \cdot\}$ denotes the commutator of operators.

Since $R(x,y)z \circ \eta =
[ (D_yD_x - D_xD_y + D_{[x,y]})z]\circ \eta$, this term is clearly
continuous in $H^s$, as $R$, the curvature of $\nabla$
on $M$, is $C^\infty$.

That $(x,y,z)\mapsto 
[ (M_yM_x - M_xM_y + M_{[x,y]})z]\circ \eta$ is continuous in $H^s$ follows
from the above property of $M$; namely, $[x,y]\in H^{s-1}(TM)$ and
for $s>(n/2)+2$, $H^{s-2}(TM)$ forms a multiplicative algebra so that
$$|M_{[x,y]} z|_{s-1} \le c \bigl|[x,y]\bigr|_{s-1}|z|_{s-1} \le
c|x|_s|y|_s|z|_s.$$

Finally, continuity of
$(x,y,z) \mapsto
\left[ \{ \nabla_x,M_y\} z + \{M_x,\nabla_y\}z\right] \circ \eta$ in $H^s$ follows
from the fact that the commutator terms are both order-zero differential operators,
together with the property of the multiplicative algebra.

\section{Smoothness of differential bundle maps over the identity}
\label{app2}

Let ${\mathfrak G}^s$ denote either ${\mathcal D}_D^s$,
${\mathcal D}_N^s$, or ${\mathcal D}_{mix}^s$.
Suppose  $L:H^s(E)\rightarrow H^{s-l}(F)$ is an order $l$ differential operator
between sections of two vector bundles $E$ and $F$ over $M$.  The purpose of 
this appendix is to carefully explain why $R_\eta \circ L \circ R_{\eta^{-1}}: 
H^s(M,E) \downarrow {\mathfrak G}^s \rightarrow H^{s-l}(M,F)\downarrow
{\mathfrak G}^s$ is smooth, even though the map $\eta \mapsto \eta^{-1}: 
{\mathfrak G}^s \longrightarrow {\mathfrak G}^s$ is only $C^0$.  That
$R_\eta \circ L \circ R_{\eta^{-1}}$ is $C^\infty$ follows from the special
structure of exact sequences covering the identity map.

A sequence of vector bundle maps over the identity
$E \stackrel{f}{\rightarrow} F \stackrel{g}{\rightarrow} G$ is {\it exact} at
$F$ if $\text{range}(f)=\text{ker}(g)$; {\it split fiber exact} if 
$\text{ker}(f)$, $\text{range}(f)$$=$$\text{ker}(g)$, and $\text{range}(g)$ 
split in $E,F$, and $G$, respectively; 
and {\it bundle exact} if additionally $\text{ker}(f)$, 
$\text{range}(f)$$=$$\text{ker}(g)$, and $\text{range}(g)$ are subbundles.
It is standard (\cite{AMR}, Proposition 3.4.20) that a split fiber exact
sequence is bundle exact, so that if $E$, $F$, and $G$ are Hilbert vector
bundles, and the sequence is exact at $F$, then $\text{ker}(f)$, 
$\text{range}(f)$$=$$\text{ker}(g)$, and $\text{range}(g)$ are subbundles.

Let $\tilde M$ denote the double of $M$, and set
$H^s(\Lambda^k) = H^s(\Lambda^k(\tilde M))$, the $H^s$ class sections of
$\Lambda^k(\tilde M)$.  Let
$H^s_\eta(\Lambda^k)$ denote the $H^s$ class maps of $\tilde M$ into 
$\Lambda^k(\tilde M)$ which cover $\eta$.

\begin{lemma} \label{smooth1}
For $s>(n/2)+1$, the map $(\eta \mapsto T\eta): {\mathcal D}^s \rightarrow
[H^s(TM)^* \otimes H^{s-1}_\eta(TM)] \downarrow {\mathcal D}^s$ is $C^\infty$.
\end{lemma}
\begin{proof}
For each $x \in M$, the metric $g$ induces a natural inner-product, say
$\bar g$, on elements of $T_x^*M \otimes T_{\eta(x)}M$, and hence a weak
$L^2$ metric on $H^s(TM)^* \otimes H^{s-1}_\eta(TM)$ given by
$\int_M \bar g (\cdot, \cdot) \mu$.  There exists a unique Levi-Civita covariant
derivative associated with this weak $L^2$ metric which we denote by
$\overline \nabla$.  The covariant derivative $\overline \nabla$ is induced by
the connector ${\mathcal K}$ which is the functorial lift of the connector
$K$ uniquely associated with the metric $\bar g$ thru the fundamental theorem
of Riemannian geometry (see Theorem 9.1 in \cite{EM}).

Let us denote the map $\eta \mapsto T\eta$ by $s$, i.e., $s(\eta) = T\eta$.
Continuity of $s$ is immediate.  Thus, we shall show that $s$ is of class
$C^1$.  Let $\epsilon \mapsto \eta^\epsilon$ be a smooth curve in 
${\mathcal D}^s$ such that $\eta^0 = \eta$ and $(d/d \epsilon)|_{\epsilon=0}
\eta^\epsilon = V_\eta \in T_\eta {\mathcal D}^s$; then,
$\overline \nabla _{V_\eta} s(\eta) \in H^s(TM)^*\otimes H^{s-1}_\eta(TM)$
is computed as
$$ \overline \nabla _{V_\eta} s(\eta) = \left.\frac{d}{d \epsilon}\right|_
{\epsilon=0} s(\eta^\epsilon) =
\left.\frac{d}{d \epsilon}\right|_{\epsilon=0} T\eta^\epsilon = \nabla V_\eta,$$
where $\nabla$ denotes the unique Levi-Civita covariant derivative in the
pull-back bundle $\eta^*(TM)$ associated to the metric $g$ on $M$. Specifically,
for $W\in T_xM$ and $V_\eta \in \eta^*(TM)$, $\nabla_W V_\eta(x)$ has the
local expression
$$ \nabla_W V_\eta(x) = TV_\eta(x) \cdot \left( T\eta(x) \cdot W(x)\right)
+ \Gamma_{\eta(x)} \left( V_\eta(x), T\eta(x) \cdot W(x)\right),$$
where $\Gamma_{\eta(x)}$ denotes the Christoffel symbol of the metric $g$
evaluated at the point $\eta(x) \in M$.

We compute the operator norm of
$\overline \nabla s(\eta) \in \text{Hom}(H^s_\eta(TM), H^s(TM)^* \otimes
H^{s-1}_\eta(TM))$ which we shall denote by $| \cdot |_{\text{op}}$.
We have that
\begin{eqnarray*}
| \overline \nabla s (\eta) |_{\text{op}} &=& \sup_{V_\eta \in H^s_\eta,
|V_\eta|_s = 1} |\nabla V_\eta|_{H^s(TM)^* \otimes H^{s-1}_\eta(TM)} \\
&=& \sup_{V_\eta \in H^s_\eta,|V_\eta|_s = 1} \sup_{W\in H^s,|W|_s = 1} 
|\nabla_W V_\eta|_{s-1} \\
&\le & \sup_{V_\eta \in H^s_\eta,|V_\eta|_s = 1} \sup_{W\in H^s,|W|_s = 1} 
|\nabla V|_{s-1} \ |W|_s \\
&<& C(g,|T\eta|_{s-1}) < \infty.
\end{eqnarray*}

Computing the supremum of $|\overline \nabla s(\eta)|_{\text{op}}$ in a 
neighborhood of $\eta$ yields the $C^1$ topology; as the supremum is finite,
we have established that $s$ is a $C^1$ map.

To see that $s$ is of class $C^2$, we compute in a local chart
\begin{align*}
\left. \frac{d}{d \epsilon}\right|_{\epsilon=0} \nabla V_{\eta^\epsilon}  =
TV_\eta (x) \cdot \nabla_W V_\eta(x) & + T \Gamma_{\eta(x)} \cdot 
T\eta(x) \left( V_\eta(x), T\eta(x) \cdot W(x)\right) \\
& + \Gamma_{\eta(x)} \left( V_\eta (x), \nabla _W V_\eta(x)\right).
\end{align*}
Since $T\eta$ is in the multiplicative algebra $H^{s-1}$, and $\Gamma \in C^
\infty$, the same argument as above shows that $s$ is $C^2$.  In particular,
we see that the $k$th derivative of $s$ is a rational combination of
$\eta, T\eta, \nabla V_\eta$ and derivatives of $\Gamma$, which combined
with our argument showing that $s$ is $C^1$ together with the fact that
multiplication of $H^{s-1}$ maps is smooth, shows that $s$ is $C^k$ for
any integer $k \ge 0$, and hence that $s$ is $C^\infty$.  
\end{proof}

Define $\overline d: H^s_\eta(\Lambda^k) \downarrow {\mathcal G}^s \rightarrow
H^{s-1}_\eta(\Lambda^{k+1}) \downarrow {\mathfrak G}^s$ to be the bundle map 
covering the identity given by
$$\overline{d}_\eta(\alpha_\eta) = [d(\alpha_\eta \circ
\eta^{-1})] \circ \eta \ \ \forall \ \alpha_\eta \in H^s_\eta(\Lambda^k).$$
Similarly, define
$\overline \delta: H^s_\eta(\Lambda^k) \downarrow {\mathcal G}^s \rightarrow
H^{s-1}_\eta(\Lambda^{k-1}) \downarrow {\mathfrak G}^s$ by
$\overline{\delta}_\eta = [\delta(\alpha_\eta \circ \eta^{-1}] \circ \eta$.
Lemma A.2 of \cite{EM} states that these bundle maps are smooth.  We give
the following proof.  First  note that, as $d$ is an antiderivation satisfying
$$ d( \alpha \wedge \beta) = d \alpha \wedge \beta + (-1)^k \alpha \wedge
d \beta \ \ \forall \ \alpha \in \Lambda^k,$$
it suffices to give the proof for $k=1$, in which case $d \alpha = \nabla \alpha
- (\nabla \alpha)^t$, where $\nabla$ is the Levi-Civita covariant derivative on
$T^*M$.  Using the chain rule, we see that $\overline d_\eta = 
[\nabla \circ T\eta^{-1} - ( \nabla \circ T\eta^{-1})^t] \circ \eta$.
Now $T\eta^{-1}$ is of class $H^{s-1}$ whenever $\eta$ is an $H^s$ class 
diffeomorphism, so the proof of Lemma \ref{smooth1} shows that $\overline d$
is $C^\infty$.  The fact that $\overline \delta$ is $C^\infty$ follows from
a similar argument.  We also have the following

\begin{lemma}\label{b1}
For $s>(n/2)+1$, if $X_\eta, Y_\eta \in H^s_\eta(T\tilde M)$, then
$$ \overline{\rm div}_\eta \circ [ \nabla (X_\eta \circ \eta^{-1}) \circ \eta
\ \cdot \  \nabla (Y_\eta \circ \eta^{-1}) \circ \eta] \in H^{s-2}_\eta
(T\tilde M).$$
\end{lemma}
\begin{proof}
We identify $X_\eta,Y_\eta \in H^s_\eta(T\tilde M)$ with 
$\alpha_\eta, \beta_\eta \in H^s_\eta (\Lambda^1)$, respectively.  
It then suffices to prove that $\overline{\delta}_\eta \circ
( \overline{d}_\eta(\alpha_\eta) \cdot \overline{d}_\eta(\beta_\eta))$ is
in $H^{s-2}_\eta(\Lambda^1)$, and hence that 
$\overline{d}_\eta(\alpha_\eta) \cdot \overline{d}_\eta(\beta_\eta)$ is
in $H^{s-1}_\eta(\Lambda^1)$ (since $\overline{\delta}$ is $C^\infty$).
But this follows since $H^{s-1}$ is a multiplicative algebra, and
$\overline{d}$ is a $C^\infty$ bundle map.
\end{proof}

A similar argument yields
\begin{lemma}\label{b1a}
For $s>(n/2)+1$, if $X_\eta, Y_\eta \in H^s_\eta(T\tilde M)$, then
$$ \overline{\rm grad}_\eta \circ \text{\rm Tr}[ \nabla (X_\eta \circ \eta^{-1}) \circ \eta
\ \cdot \  \nabla (Y_\eta \circ \eta^{-1}) \circ \eta] \in H^{s-2}_\eta
(T\tilde M).$$
\end{lemma}

We shall need Lemma A.3 in \cite{EM} which we state
as follows:

\begin{lemma} \label{lemma3}
Let $\pi: E \rightarrow M$ be a vector bundle, let ${\mathcal J}$ be
a finite dimensional subspace of $H^s(E)$ consisting of $C^\infty$ elements,
and let ${\mathfrak P}: H^s(E) \rightarrow {\mathcal J}$ be a continuous
orthogonal projector onto ${\mathcal J}$.  Then
$\overline {\mathcal J} = {\mathcal J}_\eta \downarrow {\mathcal D}^s$
is a subbundle of $H^r_\eta(M,E) \downarrow {\mathcal D}^s$ for
$r \le s$, where ${\mathcal J}_\eta = \{ f\in H^r(M,E) | f\in R_\eta 
{\mathcal J} \}$.  Furthermore, $\overline {\mathfrak P} : H^r_\eta \downarrow
{\mathcal D}^s \rightarrow \overline{\mathcal J}$, given by
$\overline{\mathfrak P}_\eta = R_\eta \circ {\mathfrak P} \circ R_{\eta^{-1}}$
is a $C^\infty$ bundle map.
\end{lemma}

For the remainder of this appendix, $\overline A$ shall denote the
bundle map given by
$\overline A_\eta(\alpha_\eta) = [A(\alpha_\eta \circ \eta^{-1})]\circ \eta$ 
for any linear operator $A$ acting on $H^s(\Lambda^k)$.
We shall use the notation $\overline{\mathcal W}$ to denote the
bundle ${\mathcal W}_\eta \downarrow {\mathcal D}^s$ for any vector space
${\mathcal W}$.  For example, $\overline{H^s(\Lambda^k)}$ shall denote
$H^s_\eta(\Lambda^k) \downarrow {\mathcal D}^s$.

Again, for $r \ge 1$, let ${\mathcal V}^r$ denote the $H^r$ vector fields
on $M$ which satisfy the boundary conditions prescribed to elements of
$T_e{\mathfrak G}^s$, and let ${\mathcal V}^r_\eta = 
\{ u \circ \eta \colon u \in {\mathcal V}^r\}$.

\begin{prop}\label{smooth_bmap}
Let ${\mathcal L} = -2\text{\rm Def}^*\text{\rm Def}$ and define 
$\overline{{\mathcal L}}$
by $\overline{{\mathcal L}}_\eta = TR_\eta \circ {\mathcal L} \circ 
TR_{\eta^{-1}}$.  Then, for $s>(n/2)+1$, and $r\ge 1$, the bundle maps
\begin{equation}\nonumber
\begin{array}{c}
\overline{(1-{\mathcal L})}:{\mathcal V}^r_\eta \downarrow {\mathfrak G}^s 
\rightarrow H^{r-2}_\eta(TM) \downarrow {\mathfrak G}^s, \\ \\
\overline{(1-{\mathcal L})^{-1}}: 
 H^{r-2}_\eta(TM) \downarrow {\mathfrak G}^s \rightarrow {\mathcal V}^r_\eta
\downarrow{\mathfrak G}^s
\end{array}
\end{equation}
are $C^\infty$.
\end{prop}

\begin{proof}
By the $L^2$ orthogonal Hodge decomposition, 
$$H^s(\Lambda^k) =
d(H^{s+1}(\Lambda^{k-1})) \oplus \delta (H^{s+1}(\Lambda^{k+1}))
\oplus {\mathcal H}^{s,k}_{\text fields},$$
where ${\mathcal H}^{s,k}_{\text fields} =\{ \alpha \in H^s(\Lambda^k) \ | \
d \alpha =0 \text{ and } \delta \alpha = 0\}$ denotes the Harmonic fields.

Hence, 
\begin{equation}\label{perp}
\left[\text{ker}(d)\right]^\perp = \delta\left(H^{s+1}(\Lambda^{k+1})\right) 
\text{ and }
\left[\text{ker}(\delta)\right]^\perp = d\left(H^{s+1}(\Lambda^{k-1})\right).
\end{equation}
Let $\pi$ denote the $L^2$ orthogonal projection of $H^{s-1}(\Lambda^{k+1})$
onto ${\mathcal H}^{s-1,k+1}_{\text fields}$, and let ${\mathfrak p} = 
\pi|_{d(H^s(\Lambda^k))}$ denote the restriction of $\pi$ to 
$d(H^s(\Lambda^k))$, so ${\mathfrak p} : d(H^s(\Lambda^k)) \rightarrow 
{\mathcal H}^{s-1,k+1}_{\text fields}$. 
Since ${\mathcal H}^{{s-1},{k+1}}_{\text fields}$ is a 
finite dimensional
subspace of $H^{s-1}(\Lambda^{k+1})$ consisting of $C^\infty$ elements, Lemma
\ref{lemma3} asserts that $\overline {\mathfrak p}$ is a smooth bundle map,
and that $\text{im}(\overline {\mathfrak p})$ and hence 
$\text{im}(1- \overline {\mathfrak p})$ is a subbundle.  We may thus form the
following exact sequence
$$
H^s_\eta(\Lambda^k) \downarrow {\mathcal D}^s 
\stackrel{\overline d}{\rightarrow} \text{im}(1-\overline {\mathfrak p})
\stackrel{\overline d}{\rightarrow} H^{s-2}_\eta(\Lambda^{k+2})\downarrow 
{\mathcal D}^s.
$$
Since $\overline d$ is a $C^\infty$ bundle map, 
this shows that $\text{ker}(\overline d)$ and $\text{im}(\overline d)$ are 
subbundles.\footnote{That $\text{ker}(d)$ and $\text{im}(d)$ are subbundles is
the statement of Lemma A.4 in \cite{EM}; we have supplied a short proof simply
to correct some typographical errors and provide some needed clarification.}

Now let ${\mathfrak p}_2: \delta(H^s(\Lambda^k)) \subset H^{s-1}(\Lambda^{k-1})
\rightarrow {\mathcal H}^{s-1,k+1}_{\text fields}$ be the restricted orthogonal
projector.
Then by the same argument $\overline {\mathfrak p}_2$ is a smooth bundle map
and $\text{im}(1- \overline {\mathfrak p}_2)$ is a subbundle.   Hence, we
may form the exact sequence
$$
H^s_\eta(\Lambda^k) \downarrow {\mathcal D}^s 
\stackrel{\overline \delta}{\rightarrow} \text{im}(1-\overline {\mathfrak p}_2)
\stackrel{\overline \delta}{\rightarrow} H^{s-2}_\eta(\Lambda^{k-2})\downarrow 
{\mathcal D}^s,
$$
and thus obtain that 
$\text{ker}(\overline \delta)$ and $\text{im}(\overline \delta)$ are subbundles.

Using (\ref{perp}), we may restrict the domain and range to ensure that
the maps
$d: \delta(H^{s+1}(\Lambda^{k+1})) \rightarrow d(H^s(\Lambda^k))$ and
$\delta: d(H^{s+1}(\Lambda^{k-1})) \rightarrow \delta(H^s(\Lambda^k))$ 
are isomorphisms.

To find the inverse of $d$ between these vector spaces, first let
$\omega = \delta \beta$.  Then
$$ d \omega = d \delta \beta \Longrightarrow \delta d \omega =
\delta d (\delta \beta) = (d \delta + \delta d) (\delta \beta) = -\triangle
\delta \beta = -\triangle \omega;$$
therefore, $\omega = (-\triangle)^{-1} \delta d \omega = \delta(-\triangle)^{-1}
d \omega,$ so that $\delta (-\triangle)^{-1}$ is the inverse of $d$.
Similarly, we find that $d(-\triangle)^{-1}$ is the inverse of $\delta$.

Next, let ${\mathfrak p}_3: \text{ker} \delta =
\delta(H^{s+1}(\Lambda^{k+1})) \oplus {\mathcal H}^{s,k}_{\text fields}
\rightarrow {\mathcal H}^{s,k}_{\text fields}$ so
$(1-{\mathfrak p}_3): \text{ker}\delta \rightarrow
\delta(H^{s+1}(\Lambda^{k+1}))$.  Now $\overline {\mathfrak p}_3$ is a smooth
bundle map by Lemma \ref{lemma3}, and since $\text{ker}(\delta)$ is a subbundle,
we may form the exact sequence
$$ \overline{\text{ker}(\delta)} 
\stackrel{\overline{\mathfrak p}_3}{\rightarrow}
\overline{{\mathcal H}^{s,k}_{\text{fields}}} \stackrel{0}{\rightarrow} 0.$$

Thus, the $\text{im}(\overline{\mathfrak p}_3)$ is a subbundle from which it
follows that $\text{im}(1-\overline{\mathfrak p}_3) = 
\overline{\delta(H^{s+1}(\Lambda^{k+1})})$ is a subbundle, so that it makes
sense to define
$$\overline d : \overline{\delta(H^{s+1}(\Lambda^{k+1}))}
\rightarrow \text{im} (\overline d)$$
as a smooth bundle isomorphism.   A similar argument allows us to define
$$\overline \delta : \overline{d(H^{s+1}(\Lambda^{k-1}))}
\rightarrow \text{im} (\overline \delta)$$ as smooth bundle isomorphism.

We have shown that the bundle map $\overline{\delta (-\triangle)^{-1}}$ 
covering the identity is the inverse of $\overline d$ which is smooth; hence,
by the inverse function theorem, the bundle map
$\overline{\delta (-\triangle)^{-1}}$  is also smooth.  On the other hand,
$\overline{d (-\triangle)^{-1}}$ is the inverse of $\overline \delta$,
and by the same argument is smooth.  Since $\overline d$ and $\overline \delta$
are $C^\infty$, then $\overline{(-\triangle)^{-1}}$ is $C^\infty$ on
$\text{im}(\overline d) \oplus  \text{im}(\overline \delta)$, and hence
$\overline{-\triangle}$ is $C^\infty$ on 
$\overline{ {\mathcal H}^{s,k}_{\text fields}}^\perp$ again by the inverse 
function theorem.

Thus far, we have been working with sections of differential $k$-forms
over the boundaryless manifold $\tilde M$.  We shall now restrict our
attention to $H^s$ class sections of $\Lambda^1(\overline M)$. Letting
$n$ denote the outward-pointing normal vector field on ${\partial M}$,
for $r\ge 2$, we define the closed subspace of $H^r(\Lambda^1(M))$ by
\begin{align*}
H^r_A = \{ \alpha \in H^r(\Lambda(M))\ |\ n \intprod &\alpha =0,
(\nabla_n \alpha^\flat)^{\rm tan} + S_n(\alpha^\flat) =0 \text{ on } 
 \Gamma_2, \\
& \text{and } \alpha =0 \text{ on } \Gamma_1 \},
\end{align*}
and for $2> r \ge 1$, set
\begin{align*}
H^r_A = \{ \alpha \in H^r(\Lambda(M))\ |\ n \intprod &\alpha =0
\text{ on } \Gamma_2, \text{ and } \alpha =0 \text{ on } \Gamma_1 \}.
\end{align*}
Note that the restriction operator to these subspaces is a continuous 
linear map.
${\mathcal L}$ is a self-adjoint linear unbounded nonnegative operator on 
$L^2$ with $D({\mathcal L}) = H^2_A$, and 
${\mathcal L}: H^2_A \rightarrow \text{im}(d) \oplus  \text{im}(\delta)$
is an isomorphism.  It follows that $(1-{\mathcal L}): H^2_A \rightarrow 
H^1(\Lambda(M))$ is an isomorphism.  Since
$${\mathcal L} = -(\triangle + 2\text{Ric} + d \delta),$$
and since we have proven that $\overline\triangle_\eta$, $\overline{d}_\eta$, 
$\overline\delta_\eta$, and $\overline{\text{Ric}}_\eta$ are $C^\infty$
bundle maps, it follows that 

$$\overline{(1-{\mathcal L})}:{(H^r_A)}_\eta \downarrow {\mathfrak G}^s 
\rightarrow H^{r-2}_\eta(TM) \downarrow {\mathfrak G}^s $$
is a $C^\infty$ bundle isomorphism covering the identity, so that
by the inverse function theorem,
$$\overline{(1-{\mathcal L})^{-1}}: 
 H^{r-2}_\eta(TM) \downarrow {\mathfrak G}^s \rightarrow {(H^r_A)}_\eta
\downarrow{\mathfrak G}^s$$
is $C^\infty$ as well.

This proves the theorem in the case that ${\mathfrak G}^s = 
{\mathcal D}_{mix}^s$.  In the case that 
${\mathfrak G}^s = {\mathcal D}_{N}^s$, simply set $\Gamma_1=\emptyset$,
and for ${\mathfrak G}^s = {\mathcal D}_{D}^s$, set $\Gamma_2=\emptyset$
in the definition of $H^r_A$.
\end{proof}

\section{Other Models of Fluid Motion}
\subsection{Third-grade fluid equations}
\label{sec9}
Set $A=\pounds_ug$ and $\alpha_1=\alpha^2$.  Let $\alpha_2 \ge 0$ and
$\beta >0$ be positive constants. The equations of a third-grade incompressible
fluid on a compact Riemannian manifold with boundary are given by
\begin{align}
\partial_t&(1-\alpha_1 \triangle_r)u - \nu \triangle_r u +
\nabla_u(1-\alpha_1\triangle_r)u - \alpha_1 (\nabla u)^t\cdot \triangle_r u
\nonumber \\
& - (\alpha_1+\alpha_2) \left(A \cdot \triangle_r u + 2 \text{div} 
\nabla u \cdot \nabla u^t\right)
-\beta \text{div}\left[\text{Tr}(A \cdot A^t)A \right] = -\text{grad } p,
\label{grade3}
\end{align}
together with the incompressibility condition $\text{div }u=0$, the
Dirichlet boundary condition $u=0$ on ${\partial M}$, and initial data
$u(0)=u_0$.
This system of equations was derived (for bounded subsets of ${\mathbb R}^n$)
by Rivlin and Ericksen \cite{RE}; equation (\ref{grade3}) generalizes the theory
to Riemannian manifolds.

For the purpose of proving well-posedness, we set all of the constants equal
to one.  It is then clear that the third-grade equations differ from equation
(\ref{avg_euler}) by the terms
$ A\triangle_r u+\text{div}\left[\nabla u \cdot \nabla u^t - \text{Tr}(A\cdot
A^t)A \right]$.
We can once again transfer the complicated study of the initial-boundary value
problem for (\ref{grade3}) to the problem of studying an ordinary differential
equation on $T{\mathcal D}_{\mu,D}^s$.  The problem of well-posedness for
this system of equations in Euclidean space has been studied previously 
in \cite{AC} and \cite{BL}.

\begin{prop}\label{3_smooth}
For $s>(n/2)+2$, let $\eta(t)$ be a curve in ${\mathcal D}_{\mu,D}^s$,
and set $u(t)=\dot\eta \circ
\eta(t)^{-1}$.  Then $u$ is a solution of the initial-boundary value problem
(\ref{grade3}) with Dirichlet boundary conditions $u=0$ on ${\partial M}$ 
if and only if
\begin{equation}\label{3_simple}
\overline{\mathcal P}_\eta \circ \left[
\frac{\nabla \dot \eta}{dt}+\left[-\nu (1-{\mathcal L})^{-1} \triangle_r u 
+ {\mathcal T}(u) + {\mathcal U}(u) + {\mathcal R}(u)\right]
\circ \eta \right] =0,
\end{equation}
where ${\mathcal U}$ and ${\mathcal R}$ are defined in Proposition
\ref{geodesic}, 
$${\mathcal T}(u) = (1-{\mathcal L})^{-1} \left[
 A\triangle_r u+\text{{\rm div}}\left(\nabla u \cdot \nabla u^t - 
\text{{\rm Tr}}(A\cdot
A^t)A \right) \right],
$$
and $\overline{\mathcal P}_\eta:T_\eta{\mathcal D}^s_D 
\rightarrow T_\eta{\mathcal D}^s_{D,\mu}$ is the Stokes projector.
\end{prop}
\begin{proof}
The proof follows trivially from Propositions \ref{geodesic}.
\end{proof}

\begin{thm}\label{3_spray}
For $s>(n/2)+2$, and $u_0 \in T_e {\mathcal D}_{\mu,D}^s $,  there exists 
$T>0$ depending on $|u_0|_s$ and independent of $\nu$, and
a unique curve $\dot \eta$ in $T{\mathcal D}^s_{D,\mu}$ satisfying
(\ref{3_simple}) with
$\eta(0) = e$ and $\dot \eta(0)=u_0$ such that
$$\dot \eta \in C^\infty( [0,T), T{\mathcal D}^s_{\mu,D})$$
has $C^\infty$  dependence on $u_0$.

For $r\ge 1$, let ${\mathcal V}^r_\mu =\{ u \in H^s(TM) \cap H^1_0(TM)\ | \
\text{\rm div }u=0\}$.
Then $u = \dot\eta \circ \eta^{-1}$  is a unique solution of the initial
value problem (\ref{grade3}), and
$$ u \in C^0( [0,T), {\mathcal V}^s_\mu) \cap C^1([0,T),{\mathcal V}^{s-1}_\mu)$$
has $C^0$ dependence on $u_0$.
\end{thm}
\begin{proof}
From Proposition \ref{3_smooth}, it is clear that
the proof is identical to the proof of Theorem \ref{Thm2} once we show that
$\dot\eta\mapsto  [{\mathcal T}(\dot\eta \circ \eta^{-1}) \circ \eta:
T_\eta{\mathcal D}_{\mu,D}^s \rightarrow T_\eta{\mathcal D}_{\mu,D}^s$
is a $C^\infty$ bundle map.   The result follows from the fact that for
$s> (n/2)+2$, $H^{s-2}$ is a multiplicative algebra, so that the terms
$A\cdot \triangle_r u$ and $\text{div}(\text{Tr}(A\cdot A^t)A)$ are of class 
$H^{s-2}$ whenever $u \in H^s$.  This observation together with the
results of Section \ref{app2} complete the proof.
\end{proof}

\subsection{A shallow water equation}
\label{sec7}

For $s>3/2$ the set
${\mathcal D}^s([0,1])$ is the Hilbert group of Dirichlet diffeomorphisms,
and $T_e{\mathcal D}^s([0,1]) = H^s(0,1)\cap H^1_0(0,1)$.

Consider the right invariant metric $\langle \cdot, \cdot \rangle$ on 
${\mathcal D}^s([0,1])$, given at the identity $e$ by
$$ \langle X, Y \rangle_e = \int_0^1 \bigl( X(x)Y(x) + X_x(x) Y_x(x)\bigl)dx.$$
As computed in \cite{M} for the group ${\mathcal D}^s(S^1)$, 
formal application of the Euler-Poincar\'{e} Theorem
\ref{thm_ep} shows that if $u(t)=\dot\eta(t) \circ \eta(t)^{-1}$, then
$\dot\eta$ is a geodesic of $\langle \cdot, \cdot\rangle$ on 
${\mathcal D}^s([0,1])$ if and only if $u(t)$ is a solution of
\begin{equation}\label{ch}
\begin{array}{c}
u_t - u_{txx} + 3uu_x -2u_xu_{xx} - uu_{xxx} =0, \\ 
u(0)=0, u(1)=0,\\
u(0)=u_0.
\end{array}
\end{equation}
This equation was derived in \cite{CH} (see also \cite{FF}).
In \cite{S}, we proved local well-posedness for the PDE (\ref{ch}) in the case
that periodic boundary conditions are imposed for all initial data $u_0$
in $H^s(S^1)$, $s>3/2$.  Our method relied on proving that the
geodesic spray of the metric $\langle \cdot, \cdot\rangle$ on 
${\mathcal D}^s(S^1)$ is smooth.  We may do the same same for on
${\mathcal D}^s([0,1])$.

\begin{thm}\label{1Dspray}
For $s>3/2$, and $u_0 \in H^s(0,1)\cap H^1_0(0,1)$,  there exists an
open interval $I=(-T,T)$, depending on $|u_0|_s$, and
a unique geodesic $\dot \eta$ of $\langle \cdot, \cdot \rangle$ satisfying
the ordinary differential equation
\begin{align*}
\ddot \eta &= {\mathcal B}(\eta,\dot\eta)=
-\bigl[(1-\partial_x^2)^{-1}\partial_x (u^2 + 
u_x^2/2)\bigr] \circ \eta,\\
\eta&(0)=e,\\
\dot\eta&(0)=u_0,\\
\end{align*}
such that
$$(\eta,\dot \eta) \in C^\infty( I,{\mathcal D}^s([0,1]) \times H^s(0,1)
\cap H^1_0(0,1))$$
has $C^\infty$  dependence on $u_0$.  

Furthermore,
$u = \dot\eta \circ \eta^{-1}$  is a unique solution of the initial
value problem (\ref{ch}), and
\begin{align*}
 u &\in C( I, H^s(0,1)\cap H^1_0(0,1)) \cap C^1(I,H^{s-1}(0,1)\cap H^1_0(0,1))
\text{ \rm if } s\ge 2, \\
 u &\in C( I, H^s(0,1)\cap H^1_0(0,1)) \cap C^1(I,H^{s-1}(0,1))
\qquad \qquad \text{ \rm if } 2>s> 3/2, \\
\end{align*}
and has $C^0$ dependence on $u_0$.
\end{thm}
\begin{proof}
For $s>3/2$, $T{\mathcal D}^s([0,1]) = {\mathcal D}^s([0,1]) \times
H^s(0,1) \cap H^1_0(0,1)$; thus to prove that $\dot\eta$ is a smooth
curve in $T{\mathcal D}^s([0,1])$, we need only copy the proof of Theorem
\ref{Thm2}, and show that ${\mathcal B}$ is a smooth map into the
second tangent bundle $T^2{\mathcal D}^s([0,1])$.  We leave the trivial
details to the reader.
\end{proof}

Having smoothness of the geodesic spray allows us to define the
Levi-Civita covariant derivative associated to $\langle\cdot, \cdot
\rangle$.

\begin{prop}\label{1Dcd}
Extending $X_\eta, Y_\eta \in T_\eta {\mathcal D}^s([0,1])$ to smooth vector 
fields
$X,Y$ on ${\mathcal D}^s([0,1])$, there exists a right invariant unique 
Levi-Civita covariant 
derivative $\tilde \nabla$ of $\langle \cdot, \cdot\rangle$ on 
${\mathcal D}^s([0,1])$
given by
\begin{align*}
\tilde \nabla_XY(\eta)& = \Bigl\{ 
\partial_t(Y_\eta\circ \eta^{-1}) + 
       \partial_x(Y_\eta\circ\eta^{-1})\cdot(X_\eta\circ \eta^{-1})\\
&\qquad \qquad
+ {\mathfrak U}(X_\eta \circ \eta^{-1}, Y_\eta \circ \eta^{-1}) 
\Bigr\} \circ \eta,
\end{align*}
where for all $u,v \in H^s(0,1)\cap H^1_0(0,1)$,
$${\mathfrak U}(u,v)=(1-\partial_x^2)^{-1}\partial_x(uv+u_xv_x/2).$$

For right-invariant vector fields $X, Y$ on ${\mathcal D}^s([0,1])$ which
are completely determined
by there value at the identity $X_e, Y_e$, 
$$
\tilde \nabla_{X}Y(e) = \partial_x(Y_e) \cdot X_e + {\mathfrak U}(X_e, Y_e).$$
\end{prop}

Again, extending $X_\eta, Y_\eta, Z_\eta \in T_\eta{\mathcal D}^s([0,1])$ to
smooth vector fields $X,Y,Z$ on ${\mathcal D}^s([0,1])$, we define
the {\it weak} Riemannian curvature tensor $\tilde R$ of the weak $H^1$ 
invariant metric $\langle
\cdot, \cdot \rangle$ on ${\mathcal D}^s([0,1])$ to be the trilinear map
$$ \tilde R_\eta : \left[T_\eta {\mathcal D}^s([0,1]))\right]^3 \rightarrow
T_\eta {\mathcal D}^s([0,1])$$
given by
$$ \tilde R_\eta(X_\eta,Y_\eta)Z_\eta= \left( \tilde \nabla_Y \tilde \nabla_X Z
\right)_\eta - \left( \tilde \nabla_X \tilde \nabla_Y Z \right)_\eta
+ \left( \tilde \nabla_{[X,Y]} Z \right)_\eta, \ \ \eta\in
{\mathcal D}^s([0,1]).  $$

Using Milnor's Lie-theoretic formula for the sectional curvature at the
identity of an invariant metric on a Lie group, Misio\l ek \cite{M} formally
computed the sectional curvature of $\tilde \nabla$ at the
identity; however the problem of showing that the weak curvature operator
$\tilde R$ is bounded in the strong $H^s$ topology was left open.  We now
establish this result.

\begin{thm}\label{1DR}
The weak curvature operator $\tilde R$ of the covariant derivative $\tilde 
\nabla$ on ${\mathcal D}^s([0,1])$ is right invariant and continuous in
the $H^s$ topology for $s>(n/2)+2$.
\end{thm}
\begin{proof}
Again, right invariance of $\tilde R$ follows from the right invariance of 
$\tilde \nabla$.

Extend $X_\eta, Y_\eta, Z_\eta \in T_\eta{\mathcal D}^s([0,1])$ to smooth 
right invariant
vector fields $x^r, y^r,z^r$ on ${\mathcal D}^s([0,1])$ and let $x=x^r(e),
y=y^r(e)$, and $z=z^r(e)$.  Let $M_xy ={\mathfrak U}(x,y)$.
Then
\begin{align*}
\tilde R_\eta(X_\eta,Y_\eta)z^r_\eta
=&
\left[ (M_yM_x - M_xM_y + M_{[x,y]})z\right]\circ \eta\\
&
\qquad + \left[ \{ \nabla_x,M_y\} z + \{M_x,\nabla_y\}z\right] \circ \eta ,
\end{align*}
where $\{ \cdot, \cdot\}$ denotes the commutator of operators, and
$\nabla_x w = (\partial_xw)\cdot x$.
Since ${\mathfrak U}(x,y)$ is in $H^s$ for $x$ and $y$ in $H^s$,
the remainder of the proof follows exactly the proof of Theorem
\ref{Thm3}.
\end{proof}

As should be clear from the above proofs, all of our results in this
section also hold for the case of periodic boundary conditions.

\section*{Acknowledgments}
The author thanks  Jerry Marsden, Marcel Oliver, and Tudor Ratiu for many 
discussions on a variety of topics appearing in this manuscript.
Research  was partially supported by the NSF-KDI grant ATM-98-73133
and the Alfred P. Sloan Foundation Research Fellowship.

\appendix
\section{The Euler-Poincar\'{e} Variational Principle}
\label{app1}

The reduction of geodesic flow on ${\mathcal D}_\mu^s$ (or any of its
subgroups) onto the single fiber of $T{\mathcal D}_\mu^s$ over the
identity $e$ is an example of the Euler-Poincar\'{e} theorem (see
\cite{MR}) which we shall now state in the setting of
a topological group $G$ which is a smooth manifold and admits smooth right
translation.  For any element $\eta$ of the group, we shall denote
by $TR_\eta$ the right translation map on $TG$, so that for example, when
$G$ is either ${\mathcal D}_{\mu,D}^s$, ${\mathcal D}_{\mu,N}^s$,  or
${\mathcal D}^s_{\mu,mix}$, then $TR_{\eta^{-1}}{\dot{\eta}} := {\dot{\eta}}
\circ \eta^{-1}$.

\begin{prop}[Euler-Poincar\'{e}]\label{thm_ep}
Let $G$ be a topological group which admits smooth manifold structure with
smooth right translation, and let $L:TG \rightarrow
{\mathbb R}$ be a right invariant Lagrangian.  Let ${\mathfrak g}$ denote
the fiber $T_eG$, and let $l:{\mathfrak g}\rightarrow{\mathbb R}$ be
the restriction of $L$ to ${\mathfrak g}$.  For a curve $\eta(t)$ in
$G$, let $u(t)=TR_{\eta(t)^{-1}}{\dot{\eta}}(t)$. Then
the following are equivalent:
\begin{itemize}
\item[\bf{a}] the curve $\eta(t)$ satisfies the Euler-Lagrange equations on
$G$;
\item[\bf{b}] the curve $\eta(t)$ is an extremum of the action function
$$S(\eta) = \int L(\eta(t),{\dot{\eta}}(t)) dt,$$
for variations $\delta \eta$ with fixed endpoints;
\item[\bf{c}] the curve $u(t)$ solves the Euler-Poincar\'{e} equations
$$ \frac{d}{dt}\frac{\delta l}{\delta u} = -\text{ad}^*_{u}
\frac{\delta l}{\delta u},$$
where the coadjoint action ad$^*_u$ is defined by
$$\langle \text{ad}^*_u v, w \rangle = \langle v, [u,w]_R\rangle,$$
for $u,v,w$ in ${\mathfrak g}$, and where $\langle \cdot, \cdot\rangle$ is
the metric on ${\mathfrak g}$ and $[ \cdot, \cdot ]_R$ is the right bracket;
\item[\bf{d}] the curve $u(t)$ is an extremum of the reduced action function
$$s(u) = \int l(u(t)) dt,$$
for variations of the form
\begin{equation}\label{con_var}
\delta u = \dot w + [w, u],
\end{equation}
where $w = TR_{\eta^{-1}}\delta \eta$ vanishes at the endpoints.
\end{itemize}
\end{prop}
See Chapter 13 in \cite{MR} for a detailed development of
the theory of Lagrangian reduction as  well as a proof of the
Euler-Poincar\'{e}
theorem.

\end{document}